\documentclass[review,hidelinks,onefignum,onetabnum]{siamart220329}

\usepackage{lipsum}
\usepackage{amsfonts}
\usepackage{graphicx}
\usepackage{epstopdf}
\usepackage{algorithmic}
\usepackage{mathtools}
\usepackage{mathrsfs}
\usepackage{xcolor}
\usepackage{enumerate}
\usepackage{amssymb}
\usepackage{subfigure}
\usepackage{longtable,booktabs}
\usepackage{lineno,hyperref}
\usepackage{threeparttable}
\newsiamthm{example}{Example}[section]

\ifpdf
  \DeclareGraphicsExtensions{.eps,.pdf,.png,.jpg}
\else
  \DeclareGraphicsExtensions{.eps}
\fi

\newsiamremark{remark}{Remark}
\newsiamremark{hypothesis}{Hypothesis}
\crefname{hypothesis}{Hypothesis}{Hypotheses}
\newsiamthm{claim}{Claim}

\headers{3D Lotka-Volterra competition models with seasonal succession }{Lei Niu, Yi Wang, and Xizhuang Xie}

\title{Global dynamics of three-dimensional Lotka-Volterra competition models with seasonal succession: I. Classification of dynamics\thanks{
\textbf{Funding:} {Niu's research was supported by the National Natural Science Foundation of
China (No. 12001096); Wang's research was supported by the National Natural Science Foundation of
China (No. 11825106 and 12331006); Xie's research was supported by the Natural Science Foundation
of Fujian Province (No. 2022J01305).}}}

\author{Lei Niu\thanks{Department of Mathematics, Donghua University, Shanghai
201620, China
  (\email{lei.niu@dhu.edu.cn}).}
\and Yi Wang\thanks{School of Mathematical Sciences, University of Science and Technology of China,
Hefei, Anhui 230026, China
  (\email{wangyi@ustc.edu.cn}).}
\and Xizhuang Xie \thanks{Corresponding author. School of Mathematical Sciences, Huaqiao University,
Quanzhou, Fujian 362021, China
(xzx@hqu.edu.cn)}. }

\usepackage{amsopn}

\ifpdf
\hypersetup{
  pdftitle={Global dynamics of three-dimensional Lotka-Volterra competition models with seasonal succession: I. Classification of dynamics},
  pdfauthor={Lei Niu, Yi Wang, and Xizhuang Xie}
}
\fi

\begin{document}

\maketitle

\begin{abstract}
The current series of two papers focus on a $3$-dimensional Lotka-Volterra competition model of differential equations with seasonal succession, which exhibits that populations experience an external periodically forced environment. We are devoted to providing a delicate global dynamical description for the model. In the first part of the series, we first use a novel technique to construct an index formula for the associated Poincar\'{e} map, by which we thoroughly classify the dynamics of the model into 33 classes via the equivalence relation relative to boundary dynamics. More precisely, we show that in classes 1--18, there is no positive fixed point and every orbit tends to certain boundary fixed point. While, for classes 19--33, there exists at least one (but not necessarily unique) positive fixed point, that is, a positive harmonic time-periodic solution of the model. Among them, the dynamics is trivial in classes 19--25 and 33, provided that the positive fixed point is unique. We emphasize that, unlike the corresponding $2$-dimensional system, a major significant difference and difficulty for the analysis of the global dynamics for $3$-dimensional system is that it may not possess the uniqueness of the positive fixed point. In the forthcoming second part of the series, we shall address the issues of (non-)uniqueness of the positive fixed points for the associated Poincar\'{e} map.
\end{abstract}



\section{Introduction}\label{section:1}
Seasonal succession is a prevalent environmental feature in nature.
The overwhelming influence of this feature on communities and ecosystems is most apparent
in temperate and high-latitude systems in which seasonal variation imposes periods of active somatic
and population growth followed by periods of depressed metabolic activity, dormancy, and increased mortality.
One striking example occurs in phytoplankton
and zooplankton of the temperate lakes, where the species grow during the warmer months and
die off or go into resting stages in winter. Such phenomenon is called seasonal succession by Sommer et al. \cite{sommer1986}.
Due to the seasonal alternate, populations experience a periodic dynamical environment driven by both
internal dynamics of species interactions and external forcing. Exploring and predicting the long-term influences of such
periodic forcing on the dynamics and structure of ecosystems is a fascinating subject (see, for example, Hale and Somolinos \cite{Hale1983},
Hess \cite{hess1991}, Hirsch and Smith \cite{hirsch2005},  Hutson et al. \cite{hutson2001}, Pol{\'a}{\v{c}}ik \cite{polacik2002} and Zhao \cite{zhao2017} ),
  and may prove to be especially vital in the future as large-scale climate change threatens to alter the strength and timing of seasonality
in many natural systems (\cite{yohe2003,walther2002}).

Klausmeier proposed a novel approach in \cite{k2010}, called successional state dynamics, to modeling succession in periodically forced food webs
and applied it to the well-known Rosenzweig-McArthur predator-prey model. He found numerically that complicated dynamics (e.g. multiannual cycles and chaos) can occur in this model.
Besides, he also set up and analyzed three other types of two-species models with seasonal succession (resource competition, facilitation and flip-flop competition). Klausmeier's approach  has been further
experimentally validated in a laboratory predator-prey experiment by Steiner et al. \cite{steiner2009}, which shows that this approach allows easy and thorough exploration of how dynamics depend on the environmental
forcing regime, and uncovers unexpected phenomena.

\par Thereafter, the ecological models with seasonal succession have attracted considerable attentions (\cite{k2012,k2013,lk2011}). However, besides some numerical and approximation results, few analytic results on the possible effect of season succession are obtained.
Mathematically, the vector fields of such models are discontinuous
and periodic with respect to time. To our knowledge, the first attempt to give an analytic study on the effect of seasonal succession was made by
Hsu and Zhao \cite{Hsu-Zhao} for the $2$-species Lotka-Volterra competition model (that is, the model without species $3$) with seasonal
succession
\begin{equation}\label{seasonal-system}
 \begin{dcases}
\frac{d x_{i}}{d t}=-\mu_{i} x_{i}, \quad t\in[k \omega,  k \omega+(1-\varphi) \omega), \\
\frac{d x_{i}}{d t}=x_{i}\bigg(b_{i}-\sum_{j=1}^{3} a_{i j} x_{j}\bigg),\quad t\in [k \omega+(1-\varphi) \omega,  (k+1) \omega),
\end{dcases}
\end{equation}
where $k \in \mathbb{Z}_{+}, \varphi \in (0,1]$ and $\omega$, $\mu_{i}$, $b_{i}$, and $a_{i j}$, $i,j=1,2,3$, are all positive constants. As a matter of fact, the model \eqref{seasonal-system} is a prototypical competition model to study the effect of seasonal succession following Klausmeier's approach \cite{k2010,k2012}, where the overall period is $\omega$, and $\varphi$ stands for the switching proportion of a period between the linear system
\begin{equation}\label{linear-system}
    \frac{d x_{i}}{d t}=-\mu_{i} x_{i}, ~~i=1, 2, 3
\end{equation}
and the classical Lotka-Volterra competition model
\begin{equation}\label{LV-system}
    \frac{d x_i}{dt}=x_i\bigg(b_i- \sum_{j=1}^3 a_{ij}x_j \bigg),~~ i=1,2,3.
\end{equation}
Biologically, seasonal succession can be thought of as a series of transitions between the two systems \eqref{linear-system} and \eqref{LV-system}, where
$\varphi$ describes the proportion of the period in good season in which the species follow the Lotka-Volterra system \eqref{LV-system}, while $(1-\varphi)$ represents the proportion of the period in bad season in which the species die exponentially according to system \eqref{linear-system}. Hsu and Zhao \cite{Hsu-Zhao} established a complete classification of the dynamics for the Poincar\'e map $\mathcal{P}$ associated with the $2$-species model \eqref{seasonal-system}. Although the phenomena obtained in \cite{Hsu-Zhao} resemble the dynamic scenarios in the classical $2$-species Lotka-Volterra competition autonomous model (see \cite{Z993}), the analysis is rather nontrivial. As a matter of fact, there are many difficulties to deal with, including the uniqueness, as well as the stability, of the fixed point of $\mathcal{P}$. By appealing to Floquet theory and the theory of monotone dynamical systems, Hsu and Zhao \cite{Hsu-Zhao} accomplished their approaches. For more recent work on seasonal succession, we refer to \cite{Niu-Wang-Xie, tang2017, xiao2016, xie2021,zhang-zhao2013}.

    However, for the $3$-dimensional model \eqref{seasonal-system} with seasonal succession, little is known about the dynamics of the model except that the associated Poincar\'e map $\mathcal{P}$ of the system admits a carrying simplex (see \cite{Niu-Wang-Xie} or \cite{Mier20,Mierczynski2023,Mier18,Wang-Jiang-02}), which is an invariant one-codimensional manifold attracting all non-trivial orbits. One may also refer to \cite{baigent2019, diekmann2008carrying, hou2020, hou2021, LG,jiang2015, ruiz2011exclusion} for the existence of carrying simplex in a series of concrete classical discrete-time competitive population models.

For time-periodically forced differential equations (including the competitive system \eqref{seasonal-system} with seasonal succession), one of the major challenge is that there is no explicit expression for its associated Poincar\'e map $\mathcal{P}$ at all. This is actually a difficult situation that will not be encountered when analyzing many concrete discrete-time competitive models in the literatures. As a consequence, compared to those well established results for the classical discrete-time competitive models (see, for example,  \cite{baigent2017, Gyllenberg2019, Gyllenberg2020b, jiang2014, LG,MNR2019}), the dynamics on the carrying simplex for competitive system \eqref{seasonal-system} still remains open and is far from understood.

The current series of two papers are devoted to a delicate global dynamical description on the carrying simplex for the $3$-dimensional competitive model \eqref{seasonal-system} with seasonal succession. To the best of our knowledge, almost all previous studies on this issue (no matter
for autonomous competitive models or concrete discrete-time competitive models) have a basic prerequisite condition, that is, the {\it uniqueness} of the positive fixed point (if exists) for competitive systems (see, e.g. \cite{baigent2013, Z993} for autonomous competitive 3D Lotka-Volterra systems, and \cite{Gyllenberg2019,Gyllenberg2020b,LG} for 3D discrete-time competitive population models). Due to the uniqueness of the positive fixed point, many key tools, such as the index formula (\cite{jiang2014,LG,MNR2019}) and criteria for ruling out nontrivial dynamics (\cite{Balreira2017,Gyllenberg2020a,Niu-Ruiz-2018}), have been introduced, which play a crucial role in studying the global dynamics of such systems.

Back to system \eqref{seasonal-system} with seasonal succession, the positive fixed point of $\mathcal{P}$ corresponds to the positive harmonic solution of \eqref{seasonal-system}. Such problems are of great attention and interest for mathematicians and theoretical ecologists. Nevertheless, even for the $2$-dimensional cases in \cite{Hsu-Zhao}, it is already quite non-trivial for verifying the uniqueness of the positive fixed point. Hence, of course, it naturally becomes a very challenge problem for the $3$-dimensional case. In the second of the series \cite{nwx2023b}, we will address the problem of uniqueness, and further find examples with a single (or multiple) positive fixed point(s), respectively.

While this paper, the first of the series,
will first make an attempt to investigate the global dynamics on the carrying simplex for the $3$-dimensional (time-periodically forced) system  \eqref{seasonal-system}
with seasonal succession. More precisely, {\it without the prerequisite condition on uniqueness of positive fixed point,}
we will provide a dynamical classification for all the associated Poincar\'{e} maps $\mathcal{P}$  via the boundary dynamics that, in particular, determines the existence  (or non-existence) of positive fixed points. Moreover, our approach is independent of the uniqueness of the positive fixed point.

To be more precise, we first establish an index formula as an effective core tool for the Poincar\'{e} map $\mathcal{P}$ (Theorem \ref{simplex-index-formula}), that is,
\begin{equation}\label{index formula}
 \sum_{\theta \in \mathscr{E}_{v}} \mathrm{ind}(\mathcal{P}, \theta)+2 \sum_{\theta \in \mathscr{E}_{s}} \mathrm{ind}(\mathcal{P}, \theta)+4 \sum_{\theta \in \mathscr{E}_{p}} \mathrm{ind}(\mathcal{P}, \theta)=1,
 \end{equation}
which reveals the intrinsic connection among all fixed points, where $\mathrm{ind}(\mathcal{P}, \theta)$ stands for the index of $\mathcal{P}$ at the fixed point $\theta$,
and $\mathscr{E}_{v}$, $\mathscr{E}_{s}$, and $\mathscr{E}_{p}$ are the sets of nontrivial fixed points on the coordinate axes, coordinate planes and the positive cone, respectively. The index formula (\ref{index formula}) for the associated Poincar\'e map $\mathcal{P}$ of system \eqref{seasonal-system} can be seen as a generalization of that provided in \cite{jiang2014} for competitive mappings, of which our approach can be applied to the classical competitive population models studied in \cite{Gyllenberg2019, Gyllenberg2020b, jiang2014, LG,ruiz2011exclusion} directly, but not vice-versa.

We then define an equivalence relation relative to local stability of fixed points on the boundary of the carrying simplex for all the $3$-dimensional Poincar\'e maps of \eqref{seasonal-system}, and derive a total of $33$ stable equivalence classes in terms of inequalities on parameters of the model  (Theorem \ref{theory:2}), and as a by-product we obtain in which classes there exists a positive fixed point. The corresponding phase portraits on the carrying simplex
are as shown in Table \ref{biao0}. Specifically, it is proved that there is no positive fixed point and every trajectory converges to some fixed point on the boundary of the carrying simplex in classes 1--18, while there is at least one (but not necessarily unique) positive fixed point in classes 19--33
(Theorem \ref{prop:3} and Proposition \ref{prop:existence-of-PFP}). Each map from class $27$ has a heteroclinic cycle, i.e. a cyclic arrangement of saddle fixed points and heteroclinic connections,  which is just the boundary of the carrying simplex. We also provide the stability criteria on the heteroclinic cycle (Corollary \ref{stable_cycle}). Moreover, we numerically find that the classes 26, 27, 29, 31 can possess attracting invariant closed curves, on which all orbits are dense, which implies the existence of attracting quasiperiodic solutions.

Here, we emphasize again that our classification DOES NOT need the prerequisite condition on the uniqueness of the positive fixed point, which is crucial in the study of numerous classical competitive differential equations and mappings (see, for example, \cite{Gyllenberg2019,Gyllenberg2020b,LG,zeeman1998three}), though the uniqueness of the positive fixed point is  very important in the study of the global dynamics for the Poincar\'e map of system \eqref{seasonal-system}. In particular, if in addition the positive fixed point is unique, then we can prove that every orbit converges to some fixed point for classes 19--25 (Corollary \ref{prop:whole-dynamics}) and the positive fixed point is globally asymptotically stable for class 33 (Corollary \ref{global-stability}).

The paper is organized as follows. In Section \ref{section:2}, we introduce our notation and provide some
 preliminaries. Section \ref{section:3} is
devoted to constructing the index formula on the carrying simplex
 for the Poincar\'e map associated with system \eqref{seasonal-system}. In Section \ref{section:4}, we define the equivalence relation relative to the boundary dynamics for all Poincar\'e maps associated with system \eqref{seasonal-system}, and derive the 33 stable equivalence classes. We list the corresponding phase portraits on the carrying simplex and parameter conditions in the Appendix. The paper ends with a
discussion in Section \ref{sec:7}.

\section{Notation and preliminaries}\label{section:2} The usual nonnegative cone of $\mathbb{R}^n$ will be denoted by $\mathbb{R}^n_+:=\{x\in \mathbb{R}^n: x_i\geq 0,~i=1,\ldots,n\}$. The interior of $\mathbb{R}^n_+$ is the open cone $\dot{\mathbb{R}}^n_+:= \{x\in \mathbb{R}^n_+: x_i>0,~i=1,\ldots,n\}$ and the boundary of $\mathbb{R}^n_+$ is $\partial \mathbb{R}^n_+:=\mathbb{R}^n_+\setminus \dot{\mathbb{R}}^n_+$. For each nonempty $I\subset \{1,\ldots,n\}$, we set $H_I^+=\{x\in \mathbb{R}_+^n: x_j=0\ \mathrm{for}\  j \notin I\}$, and $\dot{H}_I^+=\{x\in H_I^+: x_i>0 \ \mathrm{for}\  i \in I\}$. In particular, $H^+_{\{i\}}$ denotes the $i$th positive coordinate axis for $i\in I$. For $x,y\in\mathbb{R}^n$, we write $x \le  y$ if $x_i \le  y_i$ for all $1\leq i\leq n$, and $x \ll  y$ if $x_i < y_i$ for all $1\leq i\leq n$. If $x \le  y$ but $x \ne  y$ we write $x <  y$. The symbol $0$  stands for both the origin of $\mathbb{R}^n$ and the real number $0$.

Let $X\subset \mathbb{R}^n$. For a map $T: X\to X$, we denote the positive trajectory
(orbit) emanating from $y\in X$ for $T$ by the set $\{T^j(y):
j\in \mathbb{Z}_+\}$, where $T^j$ denotes the $j$-fold composition of $T$ with itself: $T\circ T\circ \cdots \circ T$, $j$ times. A set $V\subset X$ is positively invariant under $T$, if $T(V)\subset V$ and invariant if $T(V)=V$. We denote by $T|_U$ the restriction of $T$ to a subset $U\subset X$. We denote by $\mathrm{Fix}(T,U)$ the set of all fixed points of $T$ in a subset $U\subset X$.

Given a $k\times k$ square matrix $A$, we write $A\geq 0$ if $A$ is a nonnegative matrix (i.e., all its entries are nonnegative) and $A>0$ if $A$ is a positive matrix (i.e., all its entries are positive).

A map $T:\mathbb{R}^n_+\to \mathbb{R}^n_+$ is competitive (or retrotone) in a subset $\mathcal{X} \subset \mathbb{R}^n_+$ if for all $x,z\in \mathcal{X}$ with $Tx < Tz$ one has that $x_i<z_i$ provided $z_i>0$.

A {\it carrying simplex} for map $T$ is a subset $\Sigma_T$ of $\mathbb{R}_+^n\setminus \{0\}$ with the  properties:
\begin{enumerate}[(P1)]
\item $\Sigma_T$ is compact, invariant and unordered;
\item $\Sigma_T$ is homeomorphic via radial projection to the $(n-1)$-dimensional standard probability simplex $\Delta^{n-1}=\{x\in \mathbb{R}_+^n:\sum_i x_i=1\}$;
\item for any $x\in \mathbb{R}_+^n\setminus \{0\}$, there exists some $z\in \Sigma_T$ so that
$\displaystyle\lim_{k\to \infty}\|T^k(x)-T^k(z)\|=0$.
\end{enumerate}
Denote the boundary of the carrying simplex $\Sigma_T$ relative to $\mathbb{R}^n_+$ by $\partial \Sigma_T=\Sigma_T\cap \partial\mathbb{R}^n_+$  and the interior of $\Sigma_T$ relative to $\mathbb{R}^n_+$ by $\dot{\Sigma}_T=\Sigma_T \setminus \partial \Sigma_T$.

\medskip

We denote by $\Psi(t, x)$ the unique solution  of model \eqref{seasonal-system} with initial point $x\in \mathbb{R}^3$. Clearly, the domain of $\Psi(\cdot, x)$ includes $[0,+\infty)$
in case $x\in \mathbb{R}_+^3$. Since the system is $\omega$-periodic, the associated Poincar\'{e} map $\mathcal{P}$ is given by
$$\mathcal{P}(x)=\Psi(\omega, x)$$
for those $x$ for which the right-hand side is defined, namely an open set $W\subset \mathbb{R}^3$ containing $\mathbb{R}_+^3$.

Let $L:\mathbb{R}^3\to \mathbb{R}^3$ be the linear map
\begin{equation}\label{linear-map}
    x\mapsto (e^{-\mu_{1}(1-\varphi) \omega} x_{1}, e^{-\mu_{2}(1-\varphi) \omega} x_{2}, e^{-\mu_{3}(1-\varphi) \omega} x_{3})
\end{equation}
with $x \in \mathbb{R}^{3}$.
We denote by $\Phi_{t}(x)$ the solution map associated with the Lotka-Volterra competitive system \eqref{LV-system}. Then, we have
$$
\mathcal{P}(x)=\Phi_{\varphi \omega}\left(L x\right), \quad x \in W.
$$

It is easy to see that each $m$-dimensional coordinate subspace of $\mathbb{R}^{3}$ is positively invariant under $\Psi$, $(1\leq m\leq 3)$. We adopt the tradition of restricting attention to the closed nonnegative cone $\mathbb{R}_{+}^{3}$.  It follows that $H_I^+$ and $\dot{H}_I^+$ are positively invariant under $\mathcal{P}$ for each nonempty $I\subset \{1,2,3\}$. We are interested in the dynamics of the discrete-time dynamical system $\left\{\mathcal{P}^{k}\right\}_{k \geq 0}$ in $\mathbb{R}^{3}_+$.

Let $A$ be the ${3\times 3}$ matrix with entries $a_{ij}$ given in the system \eqref{seasonal-system} and set $$r_{i}=b_{i} \varphi \omega-\mu_{i}(1-\varphi) \omega,~~i,j=1,2,3.$$
We assume that $r_{i}>0$ by noticing that $(\mathcal{P}_i|_{H_{\{i\}}^+})^k(x)\to 0$ as $k\to +\infty$ for all $x\in H_{\{i\}}^+$ if $r_i\leq 0$ (see \cite[Lemma 2.1]{Hsu-Zhao}).

\begin{lemma}\label{lemma-PF}
A point $\theta \in \dot{\mathbb{R}}_{+}^{3}$ is a positive fixed point of $\mathcal{P}$ if and only if $\hat{\theta}:=\int_{0}^{\varphi \omega} \Phi_{t}(L \theta) d t$ is a positive solution of the linear algebraic system $(Ax^\tau)_i=r_i$ with $r_i>0$, $i=1,2,3$.
\end{lemma}
\begin{proof}
See the proof of Hsu and Zhao \cite[Lemma 2.5]{Hsu-Zhao}, which is easily adapted to prove the lemma.
\end{proof}

\begin{proposition}\label{lemma-eigenvalue}
Let $\theta \in \mathrm{Fix}(\mathcal{P},\partial \mathbb{R}^3_+)$ and $\hat{\theta}$ be defined in Lemma \ref{lemma-PF}. Then
$$\lambda_{i}(\theta)=\exp\Big\{r_i-\sum_{j=1}^3 a_{i j} \hat{\theta}_{j}\Big\}
$$ is an eigenvalue of $D\mathcal{P}(\theta)$ for any $i\notin  \mathcal{K}_{\theta}:=\left\{k: \theta_{k}>0\right\}$.
\end{proposition}
\begin{proof}
Assume, without loss of generality, that $\theta_3=0$. Let
\begin{equation}
    g_{i}(x)=b_{i}-\sum_{j=1}^{3} a_{i j} x_{j}
\end{equation}
and $f(x)=(f_1(x),f_2(x),f_3(x))$ with
\begin{equation}\label{equ-f}
    f_i(x)=x_ig_i(x),\quad i=1,2, 3.
\end{equation}
Let $W(t,x)=D_{x} \Phi_{t}(x)$ and $U(t,x)=D f(\Phi_{t}(x))$. Then
\begin{equation}\label{equ-DPhi}
    \frac{d W(t)}{d t}=U(t,x) \cdot W(t), \quad W(0)=I.
\end{equation}
Since $\mathcal{P}(x)=\Phi_{\varphi \omega} (Lx)$, one has
\begin{equation}\label{equ-DP}
    D \mathcal{P}(\theta)=W(\varphi \omega,L \theta) \cdot D L,
\end{equation}
where $DL=\mathrm{diag}\left(e^{-\mu_{1}(1-\varphi) \omega}, e^{-\mu_{2}(1-\varphi) \omega}, e^{-\mu_{3}(1-\varphi) \omega}\right)$. Let $v(t,x):=\Phi_{t} (x)$. Then $v_3(t,L\theta)=0$ by the expression \eqref{LV-system} and $\theta_3=0$. Thus,  $$U(t,L\theta)_{(3,3)}=g_3(\Phi_{t} (L\theta))$$ and
$$U(t,L\theta)_{(3,j)}=0~ \mathrm{for~} j\neq 3.$$
It follows that
$$W(\varphi\omega,L\theta)_{(3,3)}=\exp\Big\{\int_0^{\varphi\omega} g_3(\Phi_{t} (L\theta))dt\Big\}= \exp\Big\{b_{3} \varphi \omega-\sum_{j=1}^3 a_{3 j} \hat{\theta}_{j}\Big\}$$
and
$$W(\varphi\omega,L\theta)_{(3,j)}=0~ \mathrm{for~} j\neq 3.$$
Therefore, $\lambda_3(\theta)$ is an eigenvalue of $D\mathcal{P}(\theta)$ by \eqref{equ-DP}.
\end{proof}

\begin{proposition}\label{prop-lambda}
For any $x\in \dot{ \mathbb{R}}^3_+$, 	there exists $(D\mathcal{P}(x))^{-1}$ which is a positive matrix.
Moreover, for any $\theta \in \mathrm{Fix}(\mathcal{P},\dot{ \mathbb{R}}^3_+)$, the eigenvalue of $D\mathcal{P}(\theta)$ with the smallest modulus, say $\lambda$, satisfies $0<\lambda<1$.
\end{proposition}
\begin{proof}
By the proof of \cite[Theorem 2.3]{Niu-Wang-Xie}, we know that there exists $(D\mathcal{P}(x))^{-1}$ which is a positive matrix. Then the classical Perron-Frobenius theorem implies $\lambda$ is a simple positive eigenvalue of $D\mathcal{P}(\theta)$. Due to \eqref{equ-DPhi} and Liouville's formula,
$$
\det W(\varphi \omega, L\theta) =\exp \Big\{\sum_{i=1}^3\mu_{i}(1-\varphi) \omega-\sum_{i=1}^3\int_{0}^{\varphi \omega}a_{ii} v_{i}(t, L\theta) d t\Big\},
$$
where $v(t,L\theta)=\Phi_t(L\theta)\in \dot{\mathbb{R}}_+^3$. Then
$$
\det D\mathcal{P}(\theta)=\exp \left\{-\sum_{i=1}^3\int_{0}^{\varphi \omega}a_{ii} v_{i}(t, L\theta) d t\right\}<1,
$$
which implies that $0<\lambda<1$.
\end{proof}
\begin{lemma}[\cite{Niu-Wang-Xie}]\label{carrying-simplex}
$\mathcal{P}$ admits a carrying simplex $\Sigma_\mathcal{P}$ if $r_i>0$, $i=1,2,3$.
\end{lemma}
In this paper, we denote  the set of all maps taking $\mathbb{R}_+^3$ into itself by $\mathcal{T}(\mathbb{R}_+^3)$ and the set of all the associated Poincar\'{e} maps of the Lotka-Volterra competition models \eqref{seasonal-system}  with
seasonal succession which have a carrying simplex by $\mathrm{CLVS}(3)$. In symbols:
$$
\mathrm{CLVS}(3):=\left\{\begin{array}{l}
    \mathcal{P}\in \mathcal{T}(\mathbb{R}_+^3): \mathcal{P}=\Phi_{\varphi \omega} \circ L,~\Phi{~ \text{is~the solution~map~of}~}\eqref{LV-system}  \\
    \noalign{\smallskip}
      \text{with~} 0<\varphi<1,\omega,\mu_i,b_i,r_i,a_{ij}>0, i,j\in \{1,2,3\}
\end{array}
\right\}.
$$

\section{The index formula on the carrying simplex}\label{section:3}
The aim of this section is to develop a formula on the fixed point indices for the Poincar\'e map $\mathcal{P}$ associated with the $3$-dimensional system \eqref{seasonal-system}, which is defined on an open neighborhood $W$ of $\mathbb{R}_{+}^{3}$ in $\mathbb{R}^{3}$. For the reader's convenience, we recall some known results on the fixed point index of a continuous map (see \cite{Amann1976, Granas2003} for a more detailed discussion).

Let $U \subset \mathbb{R}^{n}$ be open and $h: U \rightarrow \mathbb{R}^{n}$ be a continuous map such that $\mathrm{Fix}(h, U)$ is compact. The fixed point index of $h$ is defined by
$$
I(h, U):=\mathrm{deg}(id-h, 0, U),
$$
where $id$ is the identity map and $\mathrm{deg}(id-h, 0, U)$ is the Brouwer degree for $id-h$. The  fixed point index of $h$ at an isolated fixed point $\theta\in U$ is defined by
$$
\mathrm{ind}(h, \theta):=I(h, B_{\delta}(\theta)),
$$
where $B_{\delta}(\theta):=\{x\in \mathbb{R}^{n}:\|x-\theta\|< \delta\}$ is an open ball in $U$ such that $\mathrm{Fix}\left(h, B_{\delta}(\theta)\right)=\{\theta\}$. In particular, if $h$ is differentiable at $\theta \in \mathrm{Fix}(h, U)$ and $1$ is not an eigenvalue of $Dh(\theta)$, then
$$
\mathrm{ind}(h, \theta)=(-1)^\beta,
$$
where $\beta$ is the sum of the multiplicities of all the eigenvalues of $Dh(\theta)$ which are greater than one.
When $h$ has only finitely many fixed points in $U$, one has
$$
I(h, U)=\sum_{\theta \in \mathrm{Fix}(h, U)} \mathrm{ind}(h, \theta).
$$

Now consider the Poincar\'e map $\mathcal{P}\in \mathrm{CLVS}(3)$. We call a fixed point $x\in \mathrm{Fix}(\mathcal{P}, \mathbb{R}_{+}^{3})$ an \textit{axial} fixed point if it lies on some coordinate axis; a \textit{planar} fixed point if it lies in the interior of some coordinate plane; and a \textit{positive} fixed point if it lies in $\dot{\mathbb{R}}^3_+$. We denote the set of all nontrivial axial, planar, and positive fixed points of $\mathcal{P}$ by $\mathscr{E}_v$, $\mathscr{E}_s$, and $\mathscr{E}_p$, respectively. Since $\mathcal{P}$ has a carrying simplex $\Sigma_\mathcal{P}$, all the nontrivial fixed points in $\mathbb{R}_+^3$ lie on $\Sigma_\mathcal{P}$, that is,
$$
\mathrm{Fix}(\mathcal{P},\Sigma_\mathcal{P})=\mathscr{E}_{v}\cup \mathscr{E}_{s}\cup \mathscr{E}_{p}.
$$
We have the following index formula on the carrying simplex for $\mathcal{P}$.
\begin{theorem}[Index Formula for $\mathcal{P}$]\label{simplex-index-formula}
Let $\mathcal{P}\in \mathrm{CLVS}(3)$. Assume that every fixed point in $\mathrm{Fix}(\mathcal{P},\mathbb{R}_{+}^{3})$ is isolated and $1$ is not an eigenvalue of $D \mathcal{P}(\theta)$ for any $\theta \in \mathrm{Fix}(\mathcal{P}, \partial \mathbb{R}_{+}^{3})$. Then all the indices of the fixed points of $\mathcal{P}$ on $\Sigma_\mathcal{P}$ satisfy that
\begin{equation}\label{equ:index}
    \sum_{\theta \in \mathscr{E}_{v}} \mathrm{ind}(\mathcal{P}, \theta)+2 \sum_{\theta \in \mathscr{E}_{s}} \mathrm{ind}(\mathcal{P}, \theta)+4 \sum_{\theta \in \mathscr{E}_{p}} \mathrm{ind}(\mathcal{P}, \theta)=1.
\end{equation}
\end{theorem}

In order to prove the index formula in \eqref{equ:index}, we point out that it is more convenient to rewrite $\mathcal{P}$ as the following  Kolmogorov type mapping:
\begin{equation}\label{equ:kolmogorov}
    \mathcal{P}\left(x_{1}, x_{2}, x_{3}\right)=\left(x_{1} F_{1}(x), x_{2} F_{2}(x), x_{3} F_{3}(x)\right)
\end{equation}
by recalling that $H_I^+$ and $\dot{H}_I^+$ are positively invariant under $\mathcal{P}$ for each nonempty $I\subset \{1,2,3\}$, where
$$
F_{i}(x):= \begin{cases}\frac{\mathcal{P}_{i}(x)}{x_{i}} & \text { if } x_{i} \neq 0, \\[3pt]
\frac{\partial \mathcal{P}_{i}}{\partial x_{i}}(x) & \text { if } x_{i}=0.\end{cases}
$$
Clearly, $F_i(x)$ are continuous functions with $F_i(x)\geq 0$ for $x\in \mathbb{R}_+^3$, $i=1,2,3$. Let $x\in \partial \mathbb{R}_{+}^{3}$. Then for any $i \in \{1,2,3\}$ such that $x_{i}=0$, one has
$$
\frac{\partial \mathcal{P}_{i}}{\partial x_{i}}(x)=F_{i}(x)\  \
{\rm and} \   \  \
\frac{\partial \mathcal{P}_{i}}{\partial x_{j}}(x)=0, \quad \forall j \neq i,
$$
which implies that $F_{i}(x)$ is an eigenvalue of $D\mathcal{P}(x)$.

Now we are ready to give the proof of Theorem \ref{simplex-index-formula}.

\smallskip

\noindent {\it Proof of Theorem \ref{simplex-index-formula}.}
Motivated by the ideas in \cite{jiang2014}, we first reflect the map $\mathcal{P}|_{\mathbb{R}_+^{3}}$ according to each coordinate plane to get a map $G$ defined on $\mathbb{R}^{3}$. In fact, by \eqref{equ:kolmogorov}, $G$ can be written as
$$
G(x)=\left(x_{1} F_{1}([x]), x_{2} F_{2}([x]), x_{3} F_{3}([x])\right),
$$
where $[x]:=(|x_1|,|x_2|,|x_3|)$.
Since the carrying simplex $\Sigma_\mathcal{P}$ is a compact set, there exists a $\sigma>0$ such that
$$
\Sigma_\mathcal{P}\subset B_\sigma \cap \mathbb{R}_{+}^3,
$$
where $B_\sigma:=\{x\in \mathbb{R}^3:\|x\|\leq \sigma\}$.
Note that there exists a $\rho>\sigma$ such that $f_{i}(x)<0$ for $i \in \mathcal{K}_{\theta}$ and $f_{i}(x)=0$ for $i\notin \mathcal{K}_{\theta}$ for those $x \in \mathbb{R}_{+}^{3}$ with $\|x\| \geq \rho$, where $f$ is defined in \eqref{equ-f}. Therefore, we have
$\Phi_{t}(x) \in \mathrm{Int} B_{\rho}$ for all $x \in B_{\rho} \cap \mathbb{R}_{+}^3$ and $t>0$,
which implies that
$$
\mathcal{P}(x)=\Phi_{\phi \omega}(L x) \in \mathrm{Int} B_{\rho},~ \forall x \in B_{\rho} \cap \mathbb{R}_{+}^{3},
$$
that is,
$$
\mathcal{P}(B_{\rho} \cap \mathbb{R}_{+}^{3}) \subset B_{\rho} \cap \mathbb{R}_{+}^{3} .
$$
Since all nontrivial fixed points of $\mathcal{P}$ in $\mathbb{R}_{+}^{3}$ lie on $\Sigma_\mathcal{P}$, $\mathcal{P}$ has only finitely many fixed points in $\mathbb{R}_{+}^{3}$ and
$$
\mathrm{Fix}(\mathcal{P}, \mathbb{R}_{+}^{3}) \subset \mathrm{Int}B_{\rho}.
$$
It follows from the definition of $G$ that
$
G(B_{\rho}) \subset B_{\rho}
$
and
$G$ has only finitely many fixed points in $\mathbb{R}^{3}$ with
$$
\mathrm{Fix}(G, \mathbb{R}^{3}) \subset \mathrm{Int} B_{\rho}.
$$
Because $\mathrm{Fix}(G,\partial B_\rho)=\emptyset$ and
$
\|G(x)\|\leq \rho =\|x\|
$
for all $x\in \partial B_\rho$, one has
$$
x\neq t G(x)\quad \mathrm{for~all~} (x,t)\in \partial B_\rho\times [0,1].
$$
Let $g:B_\rho\to \mathbb{R}^3$ be the constant map $x\mapsto 0$ and consider the homotopy $$H(x,t)=tG(x)+(1-t)g(x),\quad t\in [0,1].$$
Clearly,
$H(x,t)\neq x$ for all  $(x,t)\in \partial B_\rho \times [0,1]
$, which implies that (see \cite{Granas2003})
$$\sum_{\theta \in \mathrm{Fix}(G, \mathbb{R}^{3})} \mathrm{ind}(G, \theta)=I(G, B_\rho)=I(g,B_\rho)=1.$$
Note that if $\theta \in \mathrm{Fix}(G, \mathbb{R}^{3})$, then $[\theta] \in \mathrm{Fix}(G, \mathbb{R}_{+}^{3})$, and moreover,
$$
\mathrm{ind}(G, \theta)=\mathrm {ind}(G,[\theta])
$$
by the definition of $G$. Therefore, we have
$$
\begin{aligned}
& \mathrm{ind}(G, 0)+2 \sum_{\theta \in \mathscr{E}_{v}} \mathrm{ind}(G, \theta)+4 \sum_{\theta \in \mathscr{E}_{s}} \mathrm{ind}(G, \theta)+8 \sum_{\theta \in \mathscr{E}_{p}} \mathrm{ind}(G, \theta) \\
=& \sum_{\theta \in \mathrm{Fix}(G, \mathbb{R}^{3})} \mathrm{ind}(G, \theta) \\
=& 1.
\end{aligned}
$$
Next we prove that $\mathrm{ind}(\mathcal{P}, \theta)=\mathrm{ind}(G, \theta)$ for all $\theta \in \mathrm{Fix}(\mathcal{P},\mathbb{R}^{3}_+)$. It is clear that
$$
\mathrm{ind}(\mathcal{P}, \theta)=\mathrm{ind}(G, \theta)
$$
for all $\theta \in \mathrm{Fix}(\mathcal{P}, \dot{\mathbb{R}}_{+}^{3})$, because $\mathcal{P}|_{\mathbb{R}_{+}^{3}}=G|_{\mathbb{R}_{+}^{3}}$. For $\theta \in \mathrm{Fix}(\mathcal{P}, \partial \mathbb{R}_{+}^{3})$, one has
$$\mathcal{K}_{\theta} =\left\{i: \theta_{i}>0\right\}  \subsetneqq \{1,2,3\}.$$
By \eqref{equ:kolmogorov} and the assumptions in the theorem,
$F_{i}(\theta)\neq 1$
is an eigenvalue of $D\mathcal{P}(\theta)$ for any $i \notin \mathcal{K}_{\theta}$. It is clear that $F_{i}([\theta])=F_{i}(\theta)$ because $\theta \in \mathbb{R}_{+}^{3}$, and hence $F_{i}([\theta]) \neq 1$ for all $i \notin \mathcal{K}_{\theta}$. Then there exists  a $\delta>0$ such that $\theta$ is the unique fixed point for both $\mathcal{P}$ and $G$ in the closed $\delta$-neighborhood $\overline{B_{\delta}(\theta)}$ of $\theta$ with
$$\overline{B_{\delta}(\theta)} \subset\left\{x \in \mathbb{R}^{3}: x_{i}>0, i \in \mathcal{K}_{\theta}\right\}\cap W$$
and
\begin{equation}\label{inequ:Fi}
   \left(1-F_{i}(x)\right)\left(1-F_{i}([x])\right)>0
\end{equation}
for all $x \in \overline{B_{\delta}(\theta)}$ and $i \notin \mathcal{K}_{\theta}$.
Consider the homotopy
$$
H(x, t)=t\mathcal{P}(x)+(1-t)G(x)
$$
with $t \in[0,1]$. We claim that $H(x, t) \neq x$ for each $(x,t) \in \partial B_{\delta}(\theta)\times [0,1]$, which will imply
$$\mathrm{ind}(\mathcal{P}, \theta)=\mathrm{ind}(G, \theta).$$
If not, then there exists $(\bar{x},\bar{t})\in \partial B_{\delta}(\theta)\times (0,1)$ such that $H(\bar{x}, \bar{t})=\bar{x}$ by noticing that $\mathcal{P}(x)\neq x$ and $G(x)\neq x$ for all $x \in \partial B_{\delta}(\theta)$. Since $\mathcal{P}|_{\mathbb{R}_{+}^{3}}=G|_{\mathbb{R}_{+}^{3}}$, one has $\bar{x} \notin \mathbb{R}_{+}^{3}$, which implies that there is a $j\notin \mathcal{K}_\theta$ such that $\bar{x}_{j} \neq 0$. It follows from
$$
\bar{t}\bar{x}_{j} F_{j}(\bar{x})+(1-\bar{t})\bar{x}_{j} F_{j}([\bar{x}])=\bar{x}_{j}
$$
that
$$
0<\bar{t}(1-F_{j}(\bar{x}))^{2}+(1-\bar{t})(1-F_{j}(\bar{x}))(1-F_{j}([\bar{x}]))=0,
$$
which is a contradiction by \eqref{inequ:Fi}.
Therefore,
$$
\mathrm{ind}(\mathcal{P}, 0)+2 \sum_{\theta \in \mathscr{E}_{v}} \mathrm{ind}(\mathcal{P}, \theta)+4 \sum_{\theta \in \mathscr{E}_{s}} \mathrm{ind}(\mathcal{P}, \theta)+8 \sum_{\theta \in \mathscr{E}_{p}} \mathrm{ind}(\mathcal{P}, \theta)=1.
$$
By Proposition \ref{lemma-eigenvalue}, the eigenvalues of $D\mathcal{P}(0)$ are $e^{r_{1}}, e^{r_{2}}$ and $e^{r_{3}}$, which are all greater than one, so $\mathrm{ind}(\mathcal{P}, 0)=-1$. Now the conclusion is immediate.

\begin{remark}
   {\rm It deserves to point out that the Kolmogorov  form of Poincar\'e map $\mathcal{P}$ in \eqref{equ:kolmogorov} plays a crucial role for the proof of Theorem \ref{simplex-index-formula}. The index formula in Theorem \ref{simplex-index-formula} has been proved in \cite[Theorem 3.2]{jiang2014} under the additional assumption (called total competition as in differential equation systems; see \cite{hirsch2006}) that
   \begin{equation}\label{equ:F}
       \frac{\partial F_i(x)}{\partial x_j}<0,~~\text{for~any}~ i, j=1,2,3,
   \end{equation}
  which are frequently  satisfied for concrete discrete-time competition population models, e.g., the Ricker model \cite{Gyllenberg2019,hirsch2008}, Leslie-Gower model \cite{LG,ruiz2011exclusion}, and Atkinson-Allen model \cite{Gyllenberg2020b,jiang2014}. Actually, the total competition assumption \eqref{equ:F} enables one to construct a competitive vector field so that the Poincar\'{e}-Hopf theorem can be utilized (see \cite{jiang2014,LG}). Unfortunately, in our current situation, it is more or less hopeless to verify total competition assumption \eqref{equ:F} for the Poincar\'e map $\mathcal{P}$ in \eqref{equ:kolmogorov}. As a consequence, to obtain the index formula (\ref{equ:index}), we adopt a novel technique in the proof of Theorem \ref{simplex-index-formula} by reflecting the Poincar\'e map $\mathcal{P}$ defined on $\mathbb{R}^{3}_+$ according to each coordinate plane to obtain a continuous map defined on $\mathbb{R}^{3}$, and then utilize the Brouwer degree theory to establish index formula directly. As a matter of fact, our new approach here can also be applied to the classical competitive population models studied in \cite{Gyllenberg2019, Gyllenberg2020b, jiang2014, LG,ruiz2011exclusion}, but not vice-versa.
   }
\end{remark}

\section{Dynamics of the 3-dimensional models}\label{section:4}
In this section, we analyze the dynamics in $\mathbb{R}^{3}_+$ of the map $\mathcal{P}\in\mathrm{CLVS}(3)$. Recall that each $\mathcal{P} \in \mathrm{CLVS}(3)$ admits a 2-dimensional carrying simplex $\Sigma_\mathcal{P}$ which is homeomorphic to $\Delta^{2}$. Each coordinate plane $\Pi_i:=\{x\in \mathbb{R}^3_+:x_i=0\}$, $i=1,2,3$, is positively invariant under $\mathcal{P}$,
and $\mathcal{P}|_{\Pi_i}\in\mathrm{CLVS}(2)$ is the associated Poincar\'e map of a $2$-dimensional Lotka-Volterra competition model with seasonal succession, which admits a $1$-dimensional carrying simplex. Therefore, $\partial\Sigma_\mathcal{P}$
is composed of the three $1$-dimensional carrying simplices of $\mathcal{P}|_{\Pi_i}$, $i=1,2,3$.

For the sake of discussion, we first recall the results of Hsu and Zhao in \cite{Hsu-Zhao} of the $2$-dimensional Lotka-Volterra competition model with seasonal succession
\begin{equation}\label{seasonal-system-2d}
 \begin{dcases}
\frac{d x_{i}}{d t}=-\mu_{i} x_{i}, \quad t\in[k \omega,  k \omega+(1-\varphi) \omega), \\
\frac{d x_{i}}{d t}=x_{i}\left(b_{i}-a_{i 1} x_{1}-a_{i 2} x_{2}\right),\quad t\in [k \omega+(1-\varphi) \omega,  (k+1) \omega),
\end{dcases}
\end{equation}
where $k \in \mathbb{Z}_{+}, \varphi \in (0,1]$ and $\omega$, $\mu_{i}$, $b_{i}$, and $a_{i j}$ are all positive constants such that $r_i>0$, $i=1,2$. Besides the trivial fixed point $0$, $\mathcal{P}$ admits two axial fixed points $q_{\{1\}}=(q_1,0)$ and $q_{\{2\}}=(0,q_2)$ with $q_i>0$, $i=1,2$. Moreover,  $\mathcal{P}$ has at most one positive fixed point, say $p$, if $a_{11}a_{22}-a_{12}a_{21}\neq 0$ (see \cite[Lemma 2.5]{Hsu-Zhao}). Let
$$\hat{q}_{\{i\}}=\int_{0}^{\varphi \omega} \Phi_{t}(L q_{\{i\}}) dt,\quad  i=1,2.$$
By Lemma \ref{lemma-PF},  $\hat{q}_{\{1\}}=(\frac{r_1}{a_{11}},0)$ and $\hat{q}_{\{2\}}=(0,\frac{r_2}{a_{22}})$.

Set $\gamma_{ij}:=a_{ii}r_j-a_{ji}r_i$ for $i,j=1,2$ and $i\neq j$. By \cite[Lemma 2.3]{Hsu-Zhao} and Proposition \ref{lemma-eigenvalue}, the eigenvalue $\lambda_j(q_{\{i\}})$ determines the stability of the axial fixed point $q_{\{i\}}$, $j\neq i$. Note that $$\lambda_j(q_{\{i\}})>1~(\mathrm{resp.}~ <1)\Leftrightarrow\gamma_{i j}>0~(\mathrm{resp.} <0)$$
and
$$
\gamma_{12}\gamma_{21}>0 \Rightarrow a_{11}a_{22}-a_{12}a_{21}\neq 0.
$$
The following conclusions follow from \cite{Hsu-Zhao} and \cite{Niu-Wang-Xie} immediately.

\begin{lemma}\label{lemma:2d-axis}
Consider the Poincar\'e map $\mathcal{P}$ of model \eqref{seasonal-system-2d}. If $\gamma_{i j}>0~($resp. $<0)$, then $q_{\{i\}}$ is a saddle (resp. an asymptotically stable node), and hence repels (resp. attracts) along $\Sigma_\mathcal{P}$. Moreover, $q_{\{i\}}$ is hyperbolic if and only if $\gamma_{ij} \neq 0$.
\end{lemma}

\begin{lemma}\label{thm:CPLV-2}
Consider the Poincar\'e map $\mathcal{P}$ of model \eqref{seasonal-system-2d}.
\begin{enumerate}[{\rm(a)}]
\item If $\gamma_{12}<0,\gamma_{21}>0$, then the positive fixed point $p$ does not exist and $q_{\{1\}}$ attracts all points not on the $x_2$-axis.
\item If $\gamma_{12}>0,\gamma_{21}<0$, then the positive fixed point $p$ does not exist and $q_{\{2\}}$ attracts all points not on the $x_1$-axis.
\item If $\gamma_{12},\gamma_{21}>0$, then $\mathcal{P}$ has a hyperbolic positive fixed
point $p$ attracting all points in $\dot{\mathbb{R}}_+^2$.
\item If $\gamma_{12},\gamma_{21}<0$, then $\mathcal{P}$ has a positive fixed point $p$ which is a hyperbolic saddle. Moreover, every nontrivial trajectory tends to one of the asymptotically stable nodes $q_{\{1\}}$ or $q_{\{2\}}$ or to the saddle $p$.
\end{enumerate}
\end{lemma}
\begin{remark}\label{Positive-EVs}
The two eigenvalues of $D\mathcal{P}(\theta)$ are both positive real numbers which do not equal 1 for any $\theta\in \mathrm{Fix}(\mathcal{P}, \Sigma_\mathcal{P})$ if $\gamma_{12},\gamma_{21}\neq 0$.
\end{remark}
\subsection{Classification via the boundary dynamics}\label{subsection:4.2}

Hereafter, define the plane
$$
l_i:=\{x\in\mathbb{R}^3:a_{ii}x_i+a_{ij}x_j+a_{ik}x_k=r_i, i, j, k \mathrm{~are~distinct}\}.$$
Let $\mathbb{R}^3_+\setminus l_i=U_i \cup B_i$, where $U_i$ and $B_i$ are the unbounded and bounded disjoint components of $\mathbb{R}^3_+\setminus l_i$, respectively.

First, we analyze the possible positions of all fixed points for $\mathcal{P}
\in \mathrm{CLVS}(3)$. Besides the trivial fixed point $0$, $\mathcal{P}$ has three axial fixed points $q_{\{1\}}=(q_1,0,0)$, $q_{\{2\}}=(0,q_2,0)$, $q_{\{3\}}=(0,0,q_3)$. Let
$$\hat{q}_{\{i\}}=\int_{0}^{\varphi \omega} \Phi_{t}(L q_{\{i\}}) d t,\quad i=1,2,3.$$ By Lemma \ref{lemma-PF}, $\hat{q}_{\{i\}}$ is just the intersection of $l_i$ and the $x_i$-coordinate axis, i.e.,
\begin{equation}\label{3d-axix-fp}
    \hat{q}_{\{1\}}=(\frac{r_1}{a_{11}},0,0),~ \hat{q}_{\{2\}}=(0,\frac{r_2}{a_{22}},0),~\hat{q}_{\{3\}}=(0,0,\frac{r_3}{a_{33}}).
\end{equation}
$\mathcal{P}$ may have a planar fixed point
$$v_{\{k\}}\in \dot{\Pi}_k=\{x\in \Pi_k:x_i>0,x_j>0, i, j, k \mathrm{~are~distinct}\}.$$
If $v_{\{k\}}$ exists, we define
$$
\hat{v}_{\{k\}}=\int_{0}^{\varphi \omega} \Phi_{t}(L v_{\{k\}}) d t.
$$
By Lemma \ref{lemma-PF}, $\hat{v}_{\{k\}}$ is just the intersection of $l_i$, $l_j$ and $\Pi_k$, that is,
\begin{equation}\label{3d-planar-fp}
    (A\hat{v}_{\{i\}}^\tau)_i=r_i\mathrm{~and~} x_k=0, ~i\neq k.
\end{equation}
Moreover, Lemma \ref{lemma-PF} implies that $\mathcal{P}$ might has a positive fixed point $p\in \dot{\mathbb{R}}^3_+$ in case $l_i$, $l_j$ and $l_k$ intersect in $\dot{\mathbb{R}}^3_+$.

Set $\gamma_{ij}=a_{ii}r_j-a_{ji}r_i$,  $i,j=1,2,3$, $i\neq j$.
By Proposition \ref{lemma-eigenvalue},
\begin{equation}\label{equ-gamma-eigenvalue}
  \lambda_j(q_{\{i\}})>1~(\mathrm{resp.}~ <1)\Leftrightarrow\gamma_{i j}>0~(\mathrm{resp.} <0)
\end{equation}
and
$$
\gamma_{ij}\gamma_{ji}>0 \Rightarrow a_{ii}a_{jj}-a_{ij}a_{ji}\neq 0.
$$
Let $\beta_{ij}=\frac{a_{jj}r_i-a_{ij}r_j}{a_{ii}a_{jj}-a_{ij}a_{ji}}$
 if $\gamma_{ij}\gamma_{ji}>0$, $i,j=1,2,3$ and $i\neq j$.

\begin{lemma}\label{lemma:8}
The trivial fixed point $0$ is a hyperbolic repeller for $\mathcal{P}$.
\end{lemma}
\begin{proof}
The result is obvious, since the eigenvalues of $D\mathcal{P}(0)$ are $e^{r_1},e^{r_2},e^{r_3}>1$ by Proposition \ref{lemma-eigenvalue}.
\end{proof}
\begin{lemma}\label{lemma:axis-stability}
If $\gamma_{ij}>0$~$(<0)$ then $q_{\{i\}}$ repels $($attracts$)$ along $\partial \Sigma_\mathcal{P} \cap \Pi_k$, where $i,j,k$ are distinct. Furthermore, if $\gamma_{ij},\gamma_{ik}>0~(<0)$ then the fixed point $q_{\{i\}}$ is a repeller $($an attractor$)$ on $\Sigma_\mathcal{P};$ if $\gamma_{ij}\gamma_{ik}<0$, then the fixed point $q_{\{i\}}$ is a saddle on $\Sigma_\mathcal{P};$ and $q_{\{i\}}$ is hyperbolic if and only if $\gamma_{ij}\gamma_{ik}\neq 0$.
\end{lemma}
\begin{proof}
By noticing that $\partial \Sigma_\mathcal{P} \cap \Pi_k$ is the carrying simplex of $\mathcal{P}|_{\Pi_k}$, the assertions follow from Lemma \ref{lemma:2d-axis} and Proposition \ref{lemma-eigenvalue} directly.
\end{proof}

\begin{lemma}\label{lemma:planar-stability-1}
If $\gamma_{jk}\gamma_{kj}>0$, then $\mathcal{P}$ admits a unique fixed point $v_{\{i\}}$ in the interior of the coordinate plane $\Pi_i$, where $i,j,k$ are distinct. Moreover, if $\gamma_{jk},\gamma_{kj}<0~(>0)$, then $v_{\{i\}}$ repels $($attracts$)$ along $\partial \Sigma_{\mathcal{P}}$.
\end{lemma}
\begin{proof}
Since $\Pi_i$ is positively invariant under $\mathcal{P}$ and $\mathcal{P}|_{\Pi_i} \in \mathrm{CLVS}(2)$, the conclusions follow from Lemma \ref{thm:CPLV-2} immediately.
\end{proof}

\begin{remark}
Note that $\gamma_{ij}<0~(>0)$ if and only if $\hat{q}_{\{i\}}\in U_j~(B_j)$, $j\neq i$. The dynamics of $q_{\{i\}}$ along $\partial \Sigma_\mathcal{P} \cap \Pi_k$ can be determined by the position of $\hat{q}_{\{i\}}$ relative to the line $l_j\cap \Pi_k$, where $i,j,k$ are distinct.

\end{remark}

\begin{lemma}\label{lemma:planar-stability-2}
Suppose that there is a unique planar fixed point $v_{\{i\}}\in\dot{\Pi}_i$. Then $(A\hat{v}_{\{i\}}^\tau)_i<r_i$~$(>r_i)$ implies that $v_{\{i\}}$ locally repels $($attracts$)$ in $\dot{\Sigma}_\mathcal{P}$. Moreover, $v_{\{i\}}$ is hyperbolic if and only if $(A\hat{v}_{\{i\}}^\tau)_i\neq r_i$.
\end{lemma}
\begin{proof}
 For definiteness, we assume that $v_{\{3\}}$ exists, and say $v_{\{3\}}=(v_1,v_2,0)$. Note that $v_{\{3\}}$ is the positive
fixed point of $\mathcal{P}|_{\Pi_3}$. By Proposition \ref{lemma-eigenvalue} the Jacobian matrix for $\mathcal{P}$ at $v_{\{3\}}$ can be written as
$$
    D\mathcal{P}(v_{\{3\}})=\left (
    \begin{array}{ccc}
        D\mathcal{P}|_{\Pi_3}({v_{\{3\}}}) & * \\
        0 &  e^{r_3-(A\hat{v}_{\{3\}}^\tau)_3}
    \end{array}
    \right ).
$$
Recall Lemma \ref{thm:CPLV-2}, the dynamics of the restricted system in $\Pi_3$ is
determined by the local dynamics of $q_{\{1\}}$ and $q_{\{2\}}$. So,
$$(A\hat{v}_{\{3\}}^\tau)_3< r_3~(>r_3)\Leftrightarrow e^{r_3-(A\hat{v}_{\{3\}}^\tau)_3}>1~(<1), $$
which determines that $v_{\{3\}}$ repels (attracts) along the eigendirection not in $\Pi_3$, and hence in $\dot{\Sigma}_{\mathcal{P}}$. Moreover, the two eigenvalues of $D\mathcal{P}|_{\Pi_3}(v_{\{3\}})$ are both positive real numbers which do not equal $1$ (see Remark \ref{Positive-EVs}) if $v_{\{3\}}$ is the unique fixed point in $\dot{\Pi}_3$, so $v_{\{3\}}$ is hyperbolic if and only if $(A\hat{v}_{\{3\}}^\tau)_3\neq r_3$.\end{proof}

\begin{definition}{\rm
Two maps $\mathcal{P}, \mathcal{P}^* \in \mathrm{CLVS}(3)$ are said to be {\it equivalent relative to the boundary of the carrying simplex} if there exists a permutation $\sigma$ of $\{1,2,3\}$ such that $\mathcal{P}$ has
a fixed point $q_{\{i\}}$ {\rm(}or $v_{\{k\}}${\rm)} if and only if $\mathcal{P}^*$ has a fixed
point $q^*_{\{\sigma(i)\}}$ {\rm(}or $v^*_{\{\sigma(k)\}}${\rm)}, and
further $q_{\{i\}}$ {\rm(}or $v_{\{k\}}${\rm)} has the same hyperbolicity and local dynamics on the carrying simplex as $q^*_{\{\sigma(i)\}}$ {\rm(}or $v^*_{\{\sigma(k)\}}${\rm)}.}
\end{definition}

A map $\mathcal{P}\in\mathrm{CLVS}(3)$ is said to be {\it stable relative to the boundary of the carrying simplex} if
all the fixed points on $\partial \Sigma_\mathcal{P}$ are hyperbolic. We say that an
equivalence class is {\it stable} if each map $\mathcal{P}$ in
it is stable relative to $\partial \Sigma_\mathcal{P}$.

\begin{remark}\label{remark:stable-map}
Note that
$$(A\hat{v}_{\{k\}}^\tau)_k<r_k~(>r_k) \Leftrightarrow a_{ki}\beta_{ij}+a_{kj}\beta_{ji}<r_k~(>r_k)\Leftrightarrow \hat{v}_{\{k\}}\in B_k~(U_k).$$
A map $\mathcal{P} \in \mathrm{CLVS}(3)$ is stable relative to the boundary of the carrying simplex if and only if $\gamma_{ij} \neq 0$ and $a_{ki}\beta_{ij}+a_{kj}\beta_{ji}\neq r_k$, i.e., $(A\hat{v}_{\{k\}}^\tau)_k\neq r_k$ {\rm(}if $v_{\{k\}}$ exists{\rm)}.
\end{remark}

Given $\mathcal{P} \in \mathrm{CLVS}(3)$, consider the $3$-dimensional Leslie-Gower map
\begin{equation}\label{3d-LG}
S_\mathcal{P}: \mathbb{R}^3_+ \to \mathbb{R}^3_+, ~ (S_\mathcal{P})_i(x)=\frac{(1+r_i)x_i}{1+\sum_{j=1}^3a_{ij}x_j},~ r_i, a_{ij}>0,~ i,j=1,2,3,
\end{equation}
which has the same parameters $r_i,a_{ij}>0$ as the Poincar\'e map $\mathcal{P}$ associated with the $3$-dimensional model \eqref{seasonal-system}. By \cite{LG} we know that $\mathcal{S}_\mathcal{P}$ has a carrying simplex $\Sigma_{\mathcal{S}_\mathcal{P}}$. Besides the trivial fixed point $0$ which is a hyperbolic repeller, it is easy to see that $\mathcal{S}_\mathcal{P}$ has three axial fixed points $\hat{q}_{\{i\}}$, $i=1,2,3$, given by \eqref{3d-axix-fp}.  There is a unique fixed point $\hat{v}_k\in\partial\Sigma_{\mathcal{S}_\mathcal{P}}\cap\dot{\Pi}_k$ given by \eqref{3d-planar-fp} if and only if $\gamma_{ij}\gamma_{ji}>0$. Moreover, if $\gamma_{ji},\gamma_{ij}<0~(>0)$ then $\hat{v}_{\{k\}}$ repels $($attracts$)$ along $\partial \Sigma_{\mathcal{S}_\mathcal{P}}$ and
$(A\hat{v}_{\{k\}}^\tau)_k<r_k$~$(>r_k)$ implies that $\hat{v}_{\{k\}}$ locally repels $($attracts$)$ in $\dot{\Sigma}_{\mathcal{S}_\mathcal{P}}$.

Denote by
$$
\mathrm{CLG}(3):=\bigg\{\mathcal{S} \in \mathcal{T}(\mathbb{R}_{+}^{3}): \mathcal{S}_{i}(x)=\frac{(1+r_{i}) x_{i}}{1+\sum_{j=1}^{3} a_{i j} x_{j}},~ r_{i}, a_{i j}>0, ~i, j=1,2,3\bigg\}
$$
the set of all $3$-dimensional Leslie-Gower maps.  By the above analysis, the map $\mathcal{P}\in\mathrm{CLVS}(3)$ has an axial fixed point $q_{\{i\}}$ (resp. a planar fixed point $v_{\{k\}}$) if and only if $\mathcal{S}_\mathcal{P}\in\mathrm{CLG}(3)$ has an axial fixed point $\hat{q}_{\{i\}}$ (resp.   a planar fixed point  $\hat{v}_{\{k\}}$), and moreover,  $q_{\{i\}}$ (resp. $v_{\{k\}}$) has the same hyperbolicity and local dynamics under $\mathcal{P}$ as  $\hat{q}_{\{i\}}$ (resp. $\hat{v}_{\{k\}}$) under $\mathcal{S}_\mathcal{P}$. Furthermore, $\mathcal{P}$ is stable relative to the boundary of the carrying simplex if and only if $\mathcal{S}_\mathcal{P}$ is stable relative to the boundary of the carrying simplex, and $\mathcal{P}, \mathcal{P}^* \in \mathrm{CLVS}(3)$ are equivalent relative to the boundary
of the carrying simplex if and only if $\mathcal{S}_\mathcal{P}, \mathcal{S}_{\mathcal{P}^*} \in \mathrm{CLG}(3)$ are equivalent relative to the boundary of the carrying simplex. Therefore, the classification program in \cite{LG} (see also \cite{Gyllenberg2020b})  works for
$\mathrm{CLVS}(3)$, and $\mathrm{CLVS}(3)$ has the same equivalence classes as $\mathrm{CLG}(3)$, so we have the following conclusion.

\begin{theorem}\label{theory:2}
There are a total of $33$ stable equivalence classes in $\mathrm{CLVS}(3)$.
\end{theorem}

The carrying simplices for the $33$ stable equivalence classes in $\mathrm{CLVS}(3)$ are presented in Table \ref{biao0}. The corresponding parameter conditions of a representative element for each class are also given.

\begin{remark}\label{rmk:PFP}
By noticing that the map $\mathcal{P}$ in $\mathrm{CLVS}(3)$ has the same parameters $a_{ij}$ as the Leslie-Gower map $\mathcal{S}_\mathcal{P}$ in $\mathrm{CLG}(3)$, one has $\det A<0$ for  classes 19--25 in $\mathrm{CLVS}(3)$ while $\det A>0$ for  classes 26--33 in $\mathrm{CLVS}(3)$  by \cite{LG}.
\end{remark}

\subsection{Dynamics on the carrying simplex}
We first recall two topological results on homeomorphisms defined on a topological disk $\mathscr{D}$, that is a set  homeomorphic to the closed unit disk $$\overline{\mathcal{D}_1}=\{(x_1, x_2)\in \mathbb{R}^2: x_1^2 + x_2^2\leq 1\}.$$
A homeomorphism $h:\mathscr{D}\to \mathscr{D}$ is an orientation preserving homeomorphism if it has degree one, that is,
$$\textup{deg}(h-q_{0},0,U)=1$$
where $h(p_{0})=q_{0}$ with $p_{0}\in \mathrm{Int} \mathscr{D}$ and $U\subset \mathscr{D}$ is any open neighbourhood of $p_{0}$.

\begin{lemma}[Corollary 2.1 in \cite{ruizherreranonlinearity}]\label{t1}
Let $h: \mathscr{D}\to  \mathscr{D}$ be an orientation preserving homeomorphism defined on the topological disk $\mathscr{D}$. If $h$ has only finitely many fixed points such that
$$\mathrm{Fix}(h,\mathscr{D})\subset \partial \mathscr{D}, $$
then any trajectory of $h$ converges to some fixed point.
\end{lemma}
\begin{lemma} [Theorem 2.1 in \cite{Niu-Ruiz-2018}]\label{t2}
Let $h: \mathscr{D}\to  \mathscr{D}$ be an orientation preserving homeomorphism defined on the topological disk $\mathscr{D}$. If $h$ has only finitely many fixed points such that
$$\mathrm{Fix}(h,\mathscr{D})\cap \mathrm{Int} \mathscr{D}=\{q\}$$
with $\mathrm{ind}(h,q)=-1$, then any trajectory of $h$ converges to some fixed point.
\end{lemma}

\begin{proposition}\label{trivial-dynamics}
Assume that there are only finitely many fixed points for $\mathcal{P}\in \mathrm{CLVS}(3)$. Then every trajectory  converges to some fixed point if $\mathrm{Fix}(\mathcal{P}, \dot{\mathbb{R}}^3_+)=\emptyset$ or $\mathrm{Fix}(\mathcal{P},$ $\dot{\mathbb{R}}^3_+)=\{q\}$ with $\mathrm{ind}(\mathcal{P},q)=-1$.
\end{proposition}
\begin{proof}
Note that the map $\mathcal{P}$ has a carrying simplex $\Sigma_{\mathcal{P}}$ by Lemma \ref{carrying-simplex} which is homeomorphic to the probability simplex $\Delta^2$, so $\Sigma_{\mathcal{P}}$ is a topological disk, and moreover, $\mathcal{P}|_{\Sigma_{\mathcal{P}}}$ is an orientation preserving homeomorphism from $\Sigma_{\mathcal{P}}$ onto $\Sigma_{\mathcal{P}}$ (see \cite{ruiz2011exclusion}).
If $\mathrm{Fix}(\mathcal{P}, \dot{\mathbb{R}}^3_+)=\emptyset$, then $\mathrm{Fix}(\mathcal{P}|_{\Sigma_{\mathcal{P}}}, \dot{\Sigma}_{\mathcal{P}})=\emptyset$, and hence every trajectory on $\Sigma_{\mathcal{P}}$ converges to some fixed point by Lemma \ref{t1}. If
$\mathrm{Fix}(\mathcal{P}, \dot{\mathbb{R}}^3_+)=\{q\}$ such that $\mathrm{ind}(\mathcal{P},q)=-1$, by Proposition \ref{prop-lambda} and \cite[Corollary 4.7]{MNR2019} one has
$$
\mathrm{ind}(\mathcal{P}|_{\Sigma_{\mathcal{P}}},q)=\mathrm{ind}(\mathcal{P},q)=-1.
$$
Then Lemma \ref{t2} implies that every trajectory on $\Sigma_{\mathcal{P}}$ converges to some fixed point. Now the result follows from the property (P3) of the carrying simplex.
\end{proof}

\begin{theorem}\label{prop:3}
For each map $\mathcal{P}$ from classes 1--18, there is no positive fixed point, and every nontrivial trajectory converges to some fixed point on $\partial \Sigma_{\mathcal{P}}$.
\end{theorem}
\begin{proof}
Recall that the map $\mathcal{P}$ and Leslie-Gower map $\mathcal{S}_\mathcal{P}$ which has the same parameters as $\mathcal{P}$ have the same dynamics on the boundary of the carrying simplex, so $\mathcal{S}_\mathcal{P}$ belongs to the classes 1--18 for $\mathrm{CLG}(3)$ in \cite{LG}. Therefore,  $\mathcal{S}_\mathcal{P}$ has no positive fixed point, which implies that there is no positive solution for the linear algebraic system $(Ax^\tau)_i=r_i$, $i=1,2,3$. Then there is no positive fixed point for $\mathcal{P}$ by Lemma \ref{lemma-PF}. The conclusion now follows from Proposition \ref{trivial-dynamics}.
\end{proof}
\begin{lemma} \label{lemma:fp-index}
Assume that $\mathcal{P} \in \mathrm{CLVS}(3)$ is stable relative to the boundary of the carrying simplex. Then
\begin{enumerate}[{\rm (i)}]
\item $\mathrm{ind}(\mathcal{P},q_{\{i\}})=1$ {\rm(}resp.\ $\mathrm{ind}(\mathcal{P},v_{\{k\}})=1${\rm)} if $q_{\{i\}}$ {\rm(}resp.\ $v_{\{k\}}${\rm)} is a repeller or an attractor on $\Sigma_\mathcal{P};$
\item $\mathrm{ind}(\mathcal{P},q_{\{i\}})=-1$ {\rm(}resp.\ $\mathrm{ind}(\mathcal{P},v_{\{k\}})=-1${\rm)} if $q_{\{i\}}$ {\rm(}resp.\ $v_{\{k\}}${\rm)} is a saddle on $\Sigma_\mathcal{P}$.
\end{enumerate}
\end{lemma}
\begin{proof}
It follows from Remark \ref{Positive-EVs} and Lemmas \ref{remark:stable-map} and \ref{lemma:planar-stability-2} that all the eigenvalues of $D\mathcal{P}(q_{\{i\}})$ and $D\mathcal{P}(v_{\{k\}})$ (if $v_{\{k\}}$ exists) are positive real numbers and do not equal $1$ if $\mathcal{P} \in \mathrm{CLVS}(3)$ is stable relative to the boundary of the carrying simplex. If $q_{\{i\}}$ (resp. $v_{\{k\}}$) is a repeller or an attractor on $\Sigma_\mathcal{P}$ then the number of the eigenvalues of $D\mathcal{P}(q_{\{i\}})$ (resp. $D\mathcal{P}(v_{\{k\}})$) greater than $1$ is even, and hence $\mathrm{ind}(\mathcal{P},q_{\{i\}})=1$ (resp. $\mathrm{ind}(\mathcal{P},v_{\{k\}})=1$). If $q_{\{i\}}$ (resp. $v_{\{k\}}$) is a saddle on $\Sigma_\mathcal{P}$ then the number of the eigenvalues of $D\mathcal{P}(q_{\{i\}})$ (resp. $D\mathcal{P}(v_{\{k\}})$) greater than $1$ is odd, and hence $\mathrm{ind}(\mathcal{P},q_{\{i\}})=-1$ (resp. $\mathrm{ind}(\mathcal{P},v_{\{k\}})=-1$).
\end{proof}

We set $\mathrm{ind}(\mathcal{P},v_{\{k\}})=0$ (resp. $\mathrm{ind}(\mathcal{P},p)=0$) if the planar (resp. positive) fixed point $v_{\{k\}}$ (resp. $p$) does not exist.
By using the index formula (Theorem \ref{simplex-index-formula}), we can obtain the existence of the positive fixed point for $\mathcal{P}$ in classes 19--33.
\begin{proposition}\label{prop:existence-of-PFP}
For each map $\mathcal{P}$ from classes 19--33, there is at least one positive fixed point.
\end{proposition}
\begin{proof}
The result can be obtained directly by Theorem \ref{simplex-index-formula}. Here we take the class $20$ as an example, since other cases can be proved similarly.
For map $\mathcal{P}$ in class $20$, one has
 $$\mathrm{Fix}(\mathcal{P},\partial \Sigma_\mathcal{P}) =\{q_{\{1\}},q_{\{2\}},q_{\{3\}},v_{\{1\}},v_{\{3\}}\},$$
of which $q_{\{1\}}$, $q_{\{2\}}$ are attractors,  $q_{\{3\}}$ is a saddle and $v_{\{1\}}, v_{\{3\}}$ are repellers  on $\Sigma_\mathcal{P}$; see Fig. \ref{fig:c20}.
\begin{figure}[h]
 \begin{center}
 \includegraphics[width=0.32\textwidth]{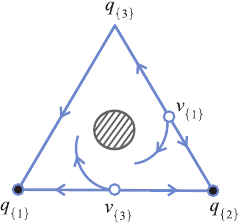}
\caption{The phase portrait on the carrying simplex for class $20$. The fixed point notation is as in Table \ref{biao0}.} \label{fig:c20}
  \end{center}
\end{figure}
\\Therefore, by Lemma \ref{lemma:fp-index} we have
$$\mathrm{ind}(\mathcal{P},q_{\{1\}})=\mathrm{ind}(\mathcal{P},q_{\{2\}})=\mathrm{ind}(\mathcal{P},v_{\{1\}})=\mathrm{ind}(\mathcal{P},v_{\{3\}})=1$$
and
$$
\mathrm{ind}(\mathcal{P},q_{\{3\}})=-1.
$$
Assume by contradiction that there is no positive fixed point for $\mathcal{P}$. Then
$$
 5=\sum^3_{i=1}(\mathrm{ind}(\mathcal{P},q_{\{i\}})+2\mathrm{ind}(\mathcal{P},v_{\{1\}})+2\mathrm{ind}(\mathcal{P},v_{\{3\}}))=1
$$
by Theorem \ref{simplex-index-formula}, which is a contradiction.
\end{proof}

Suppose that a $3$-dimensional map
$S=(x_{1} F_{1}(x), x_{2} F_{2}(x), x_{3} F_{3}(x))$
has a carrying simplex $\Sigma_S$ such that $S(\partial \Sigma_S)\subset \partial \Sigma_S$ and $S(\dot{\Sigma}_S)\subset \dot{\Sigma}_S$. Suppose further that $q_{\{1\}}=(q_1,0,0)$,
$q_{\{2\}}=(0,q_2,0)$ and $q_{\{3\}}=(0,0,q_3)$ are its three axial fixed points lying on the vertices of $\Sigma_S$ and $F_j(q_{\{i\}})>0$ for $i\neq j$. If each $q_{\{i\}}$ is a saddle, and $\partial \Sigma_S\cap \pi_i$ is the saddle connection between $q_{\{j\}}$ and $q_{\{k\}}$, then $S$ admits a heteroclinic cycle of May-Leonard type  (\cite{CHW,jiang2015,may1975}): $q_{\{1\}}\to q_{\{2\}} \to q_{\{3\}}\to q_{\{1\}}$ (or the arrows reserved), which is just the boundary of $\Sigma_S$.

Note that for any map $\mathcal{P}$ in class $27$, each axial fixed point $q_{\{i\}}$ is a saddle on $\Sigma_\mathcal{P}$, and $\partial \Sigma_\mathcal{P}\cap \Pi_i$ is the heteroclinic connection between saddles $q_{\{j\}}$ and $q_{\{k\}}$, where $i,j,k$ are distinct. So $\partial \Sigma_\mathcal{P}$ forms a heteroclinic cycle of the May-Leonard type, i.e., any map $\mathcal{P}$ in class $27$ admits a heteroclinic cycle (see Table \ref{biao0} (27)).

\begin{corollary}\label{stable_cycle}
Assume that $\mathcal{P}\in\mathrm{CLVS}(3)$ is in class $27$. If
\begin{equation}\label{equ:vartheta}
\vartheta:=w_{12} w_{23} w_{31}+w_{21}w_{13}w_{32}<0~(resp.~>0),
\end{equation}
where $w_{ij}=r_j-a_{ji}\frac{r_i}{a_{ii}}$,  $i,j=1,2,3$, $i\neq j$, then the heteroclinic cycle $\partial \Sigma_\mathcal{P}$ of $\mathcal{P}$ attracts {\rm(}resp.~repels{\rm )}.
\end{corollary}
\begin{proof}
By \eqref{equ:kolmogorov}, the map $\mathcal{P}$ can be written as
$$
\mathcal{P}(x)=(x_1F_1(x),x_2F_2(x),x_3F_3(x)),
$$
and moreover,
$$
F_j(q_{\{i\}})=\lambda_j(q_{\{i\}})=e^{r_j-a_{ji}\frac{r_i}{a_{ii}}}
$$
for $i,j=1,2,3$ and $i\neq j$ by Proposition \ref{lemma-eigenvalue}. Then the conclusion follows from \cite[Theorem 3]{jiang2015} immediately.
\end{proof}
\subsubsection{Global dynamics when the positive fixed point is unique}
It deserves to point out that by our classification in Theorem \ref{theory:2}, one can obtain the global dynamics of classes 19--25 and 33 immediately whenever they have a unique positive fixed point.  Specifically, we have the following conclusions.

\begin{corollary}\label{prop:whole-dynamics}
	If the map $\mathcal{P}$ from classes 19--25 has a unique positive fixed point, say $p$, then every trajectory converges to some fixed point. Moreover, if $1$ is not an eigenvalue of $D\mathcal{P}(p)$, then $p$ is a hyperbolic saddle whose stable manifold and unstable manifold are simple curves; in this case, the phase portraits on the carrying simplices for these classes are as shown in Fig. \ref{fig:c19-25-33}.
\end{corollary}
\begin{proof}
	By the similar arguments in the proof of Proposition \ref{prop:existence-of-PFP}, one can prove
	$$\mathrm{ind}(\mathcal{P},p)=-1$$
	if
	$\mathrm{Fix}(\mathcal{P}, \dot{\mathbb{R}}^3_+)=\{p\}$. It then follows from Proposition \ref{trivial-dynamics} that every trajectory converges to some fixed point.
	Moreover, if $1$ is not an eigenvalue of $D\mathcal{P}(p)$, then together with Proposition \ref{prop-lambda} we know that all three eigenvalues of $D\mathcal{P}(p)$, say $\lambda, \lambda_1,\lambda_2$, are positive real numbers with
	$0<\lambda<\lambda_1<1<\lambda_2$.
	Therefore, $p$ is a hyperbolic saddle whose stable manifold and unstable manifold on the carrying simplex are simple curves (see \cite{MNR2019}). Now the whole dynamics on the carrying simplices is immediate which is as shown in Fig. \ref{fig:c19-25-33} (see \cite[Corollary 5.4]{MNR2019}).
\end{proof}

\begin{corollary}\label{global-stability}
If the map $\mathcal{P}$ in class 33 has a unique positive fixed point $p$, then $p$ is globally asymptotically stable in $\dot{\mathbb{R}}^3_+$, and the phase portrait on the carrying simplex for the class 33 is as shown in Fig. \ref{fig:c19-25-33}.
\end{corollary}
\begin{proof}
By the proof of Proposition \ref{prop-lambda}, we know that $\det D\mathcal{P}(x)>0$ for all $x\in \mathbb{R}^3_+$ and $D\mathcal{P}(x)^{-1}>0$ for all $x\in \dot{\mathbb{R}}^3_+$.  By the condition (i) in Table \ref{biao0} (33), we know that $\gamma_{kj},\gamma_{jk}>0$ which implies that there exists a unique fixed point $v_{\{i\}}\in \dot{\Pi}_i$, where $i,j,k$ are distinct. Moreover, Lemma \ref{thm:CPLV-2} implies that $v_{\{i\}}$ is globally asymptotically
stable in the interior of $\Pi_i$. Since $a_{ij}\beta_{jk}+a_{ik}\beta_{kj}<r_i$, i.e., $(A\hat{v}_{\{i\}}^\tau)_i<r_i$, it follows from Lemma \ref{lemma:planar-stability-2} that $v_{\{i\}}$ locally repels in $\dot{\Sigma}_\mathcal{P}$, and hence $v_{\{i\}}$ is a saddle for $\mathcal{P}$. Then the conclusion follows from \cite[Theorem 2.4]{Balreira2017} (or \cite[Theorem 1.2]{Gyllenberg2020a}) immediately if $\mathcal{P}$ has a unique positive fixed point $p$.
\end{proof}

\subsection{Numerical simulation}\label{sec:cycle}
An interesting problem is whether there are invariant closed curves for the associated Poincar\'e map of system \eqref{seasonal-system}.
Our simulations in this section illustrate numerically that there do exist attracting invariant closed curves in some classes, such as classes 26, 27, 29, 31.

\begin{figure}[H]
    \centering
    \begin{tabular}{cc}
        \subfigure[The solution of system (\ref{seasonal-system})]{
            \label{class26-a}
            \begin{minipage}[b]{0.4\textwidth}
                \centering
                \includegraphics[width=\textwidth]{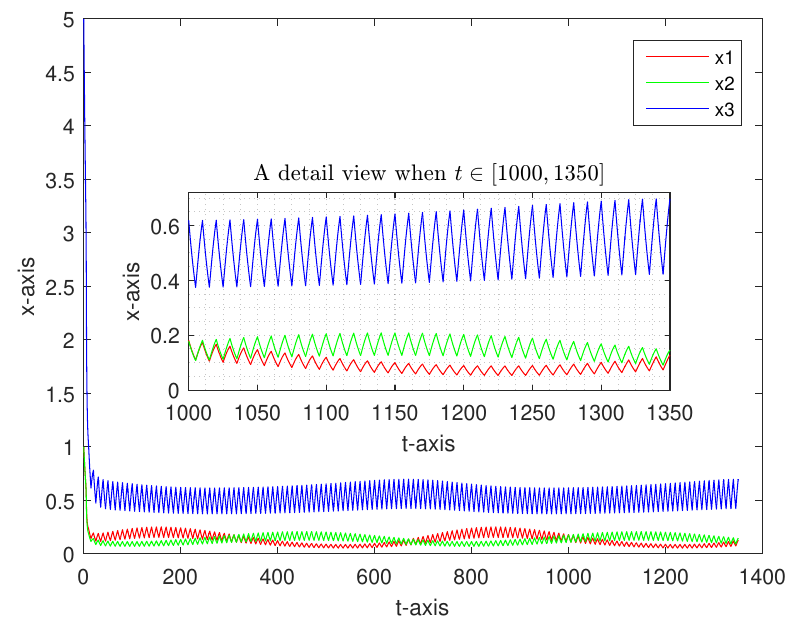}
            \end{minipage}
        } &
        \subfigure[The orbit of map $\mathcal{P}$]{
            \label{class26-b}
            \begin{minipage}[b]{0.42\textwidth}
                \centering
                \includegraphics[width=\textwidth]{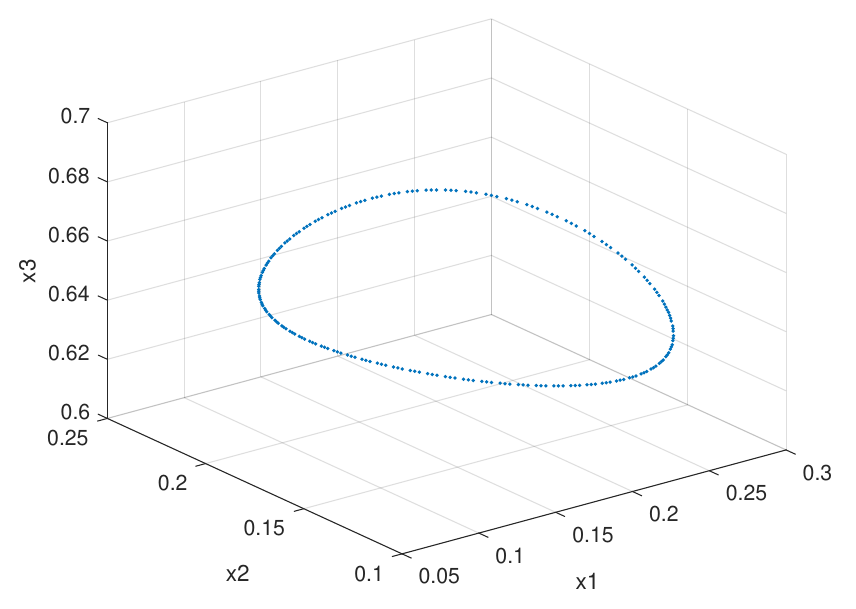}
            \end{minipage}
        } \\
    \end{tabular}
    \caption{The attracting invariant closed curve can occur in class 26.
    }
    \label{fig-c26}
\end{figure}
\begin{example}[Invariant closed curves in the class 26]
{\rm Taking parameter values
$\omega=10$, $\varphi=0.5$, $\mu_1=0.1$, $\mu_2=0.1$, $\mu_3=0.1$, $b_1=0.3$
, $b_2=0.3$, $b_3=0.3$, $a_{11}=0.3$, $a_{12}=0.6$, $a_{13}=0.15$, $a_{21}=0.1$, $a_{22}=0.2$, $a_{23}=0.3$,
$a_{31}=0.2$, $a_{32}=0.3$, $a_{33}=0.25$,
system (\ref{seasonal-system}) satisfies the inequalities of the class 26 in Table \ref{biao0}. The numerical simulations for the solution of system \eqref{seasonal-system}
with initial value $x_0=(1, 1,5)$ and the orbit of the associated Poincar\'e map $\mathcal{P}$ are shown in Fig. \ref{fig-c26}, which imply $\mathcal{P}$
admits an attracting invariant closed curve on $\Sigma_\mathcal{P}$.}
\end{example}

   \begin{figure}[H]
    \centering
    \begin{tabular}{cc}
        \subfigure[The solution of system (\ref{seasonal-system})]{
            \label{class27-2-a}
            \begin{minipage}[b]{0.4\textwidth}
                \centering                \includegraphics[width=\textwidth]{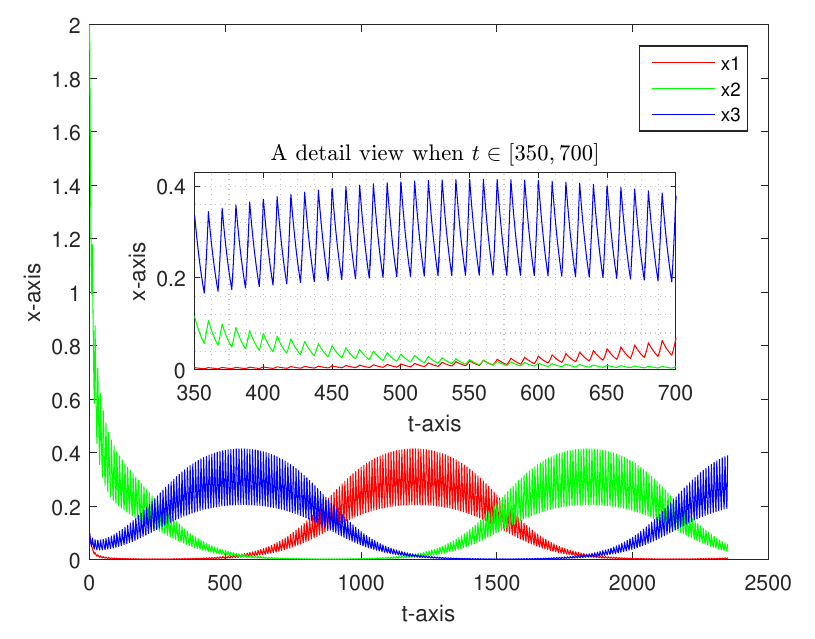}
            \end{minipage}
        } &
        \subfigure[The orbit of map $\mathcal{P}$]{
            \label{class27-2-b}
            \begin{minipage}[b]{0.42\textwidth}
                \centering                \includegraphics[width=\textwidth]{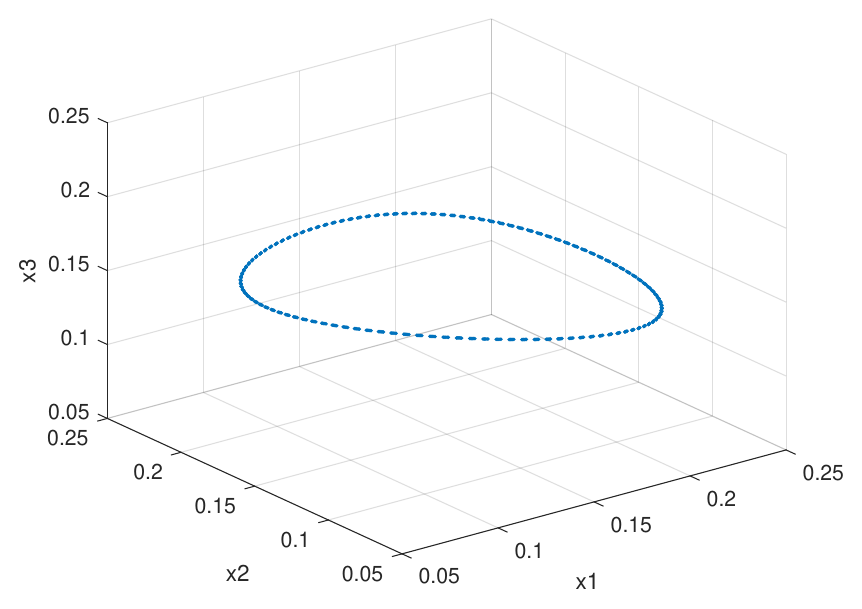}
            \end{minipage}
        } \\
    \end{tabular}
    \caption{The attracting invariant closed curve can occur in class 27.
    }
    \label{fig-c27-2}
\end{figure}
\begin{example}[Invariant closed curves in the class 27]
{\rm Taking parameter values
$\omega=10$, $\varphi=0.3$, $\mu_1=0.1$, $\mu_2=0.1$, $\mu_3=0.1$, $b_1=0.3$
, $b_2=0.3$, $b_3=0.3$, $a_{11}=0.2$, $a_{12}=0.3$, $a_{13}=0.1$, $a_{21}=0.1$, $a_{22}=0.2$, $a_{23}=0.3$,
$a_{31}=0.3$, $a_{32}=0.1$, $a_{33}=0.2$,
system (\ref{seasonal-system}) satisfies the inequalities of class 27 in Table \ref{biao0}. The numerical simulations for the solution of system \eqref{seasonal-system} with initial value $x_0=(1,2,1)$
and the orbit of the associated Poincar\'e map $\mathcal{P}$ are shown in Fig.  \ref{fig-c27-2}, which imply that the given system admits
 an attracting invariant closed curve.}
 \end{example}

 \begin{figure}[H]
    \centering
    \begin{tabular}{cc}
        \subfigure[The solution of system (\ref{seasonal-system})]{
            \label{class29-a}
            \begin{minipage}[b]{0.4\textwidth}
                \centering                \includegraphics[width=\textwidth]{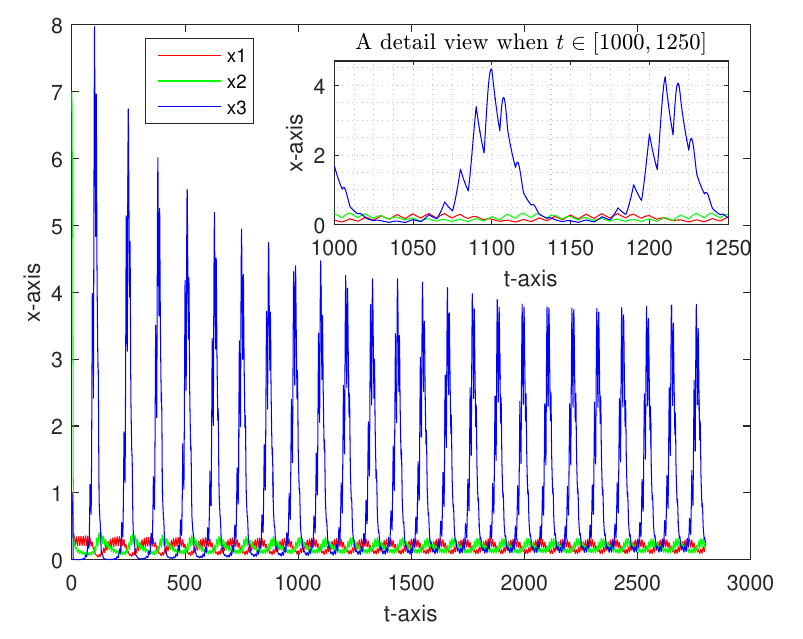}
            \end{minipage}
        } &
        \subfigure[The orbit of map $\mathcal{P}$]{
            \label{class29-b}
            \begin{minipage}[b]{0.42\textwidth}
                \centering                \includegraphics[width=\textwidth]{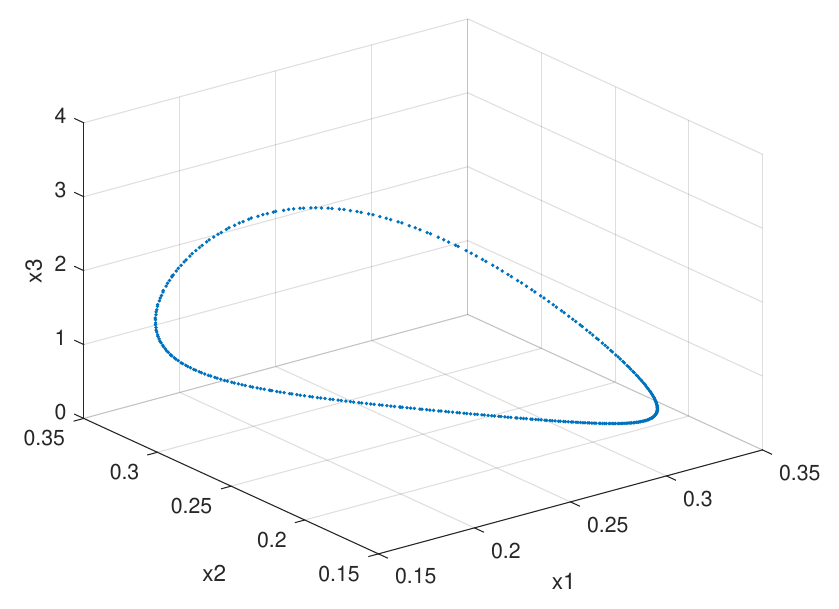}
            \end{minipage}
        } \\
    \end{tabular}
    \caption{The attracting invariant closed curve can occur in class 29.
    }
    \label{fig-c29}
\end{figure}

\begin{example}[Invariant closed curves in the class 29]
{\rm Taking parameter values
$\omega=10$, $\varphi=0.5$, $\mu_1=0.1$, $\mu_2=0.1$, $\mu_3=0.1$, $b_1=0.34$
, $b_2=0.5$, $b_3=0.7$, $a_{11}=238/325$, $a_{12}=73/325$, $a_{13}=14/325$, $a_{21}=357/325$, $a_{22}=292/325$, $a_{23}=1/325$,
$a_{31}=119/650$, $a_{32}=73/26$, $a_{33}=3/325$,
system (\ref{seasonal-system}) satisfies the inequalities of class 29 in Table \ref{biao0}. The numerical simulations for the solution of system \eqref{seasonal-system} with initial value $x_0=(1,7,1)$
and the orbit of the associated Poincar\'e map $\mathcal{P}$ are shown in Fig.  \ref{fig-c29}, which imply the given system
admits an attracting invariant closed curve.}
\end{example}
 \begin{figure}[h!]
    \centering
    \begin{tabular}{cc}
        \subfigure[The solution of system (\ref{seasonal-system})]{
            \label{class31-a}
            \begin{minipage}[b]{0.4\textwidth}
                \centering
                \includegraphics[width=\textwidth]{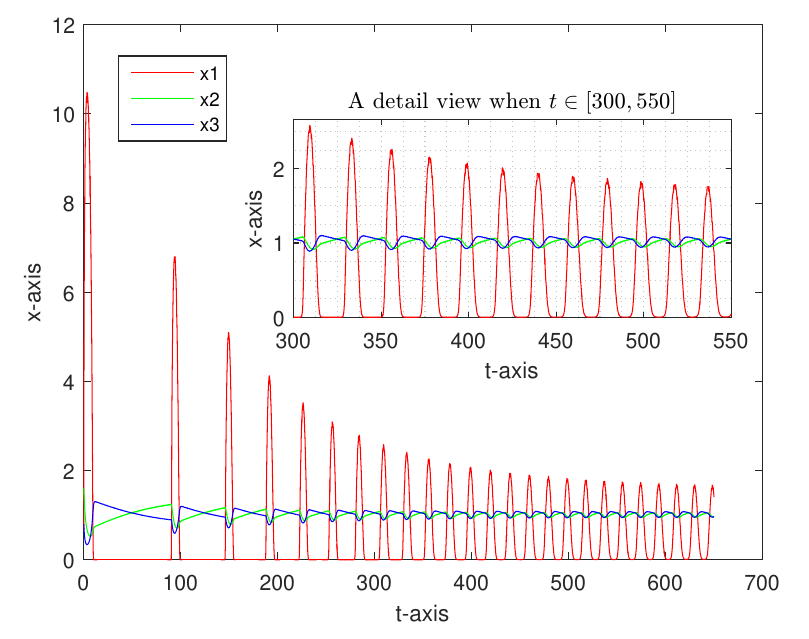}
            \end{minipage}
        } &
        \subfigure[The orbit of map $\mathcal{P}$]{
            \label{class31-b}
            \begin{minipage}[b]{0.42\textwidth}
                \centering
                \includegraphics[width=\textwidth]{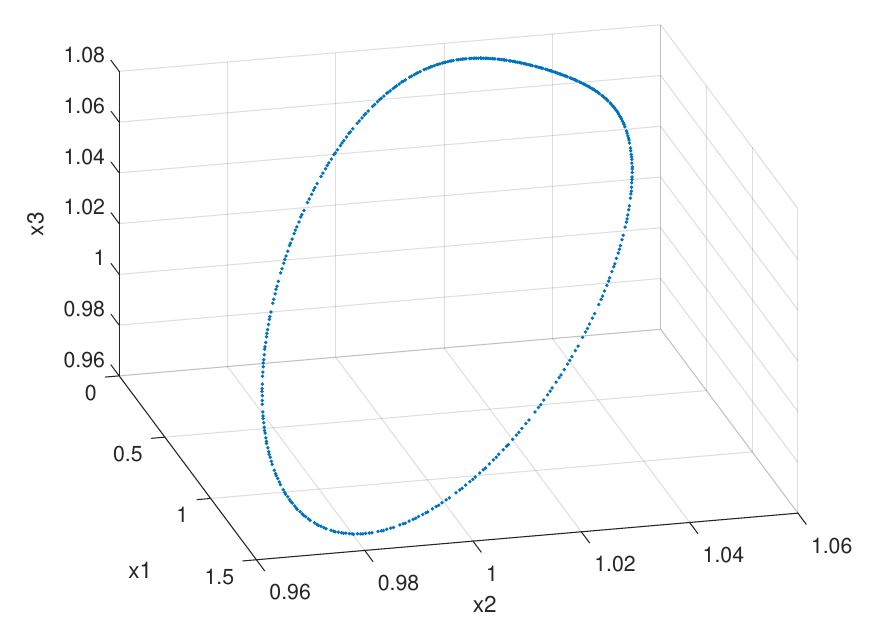}
            \end{minipage}
        } \\
    \end{tabular}
    \caption{The attracting invariant closed curve can occur in class 31.
    }
    \label{fig-c31}
\end{figure}
\begin{example}[Invariant closed curves in the class 31]
{\rm Taking parameter values
$\omega=1$, $\varphi=0.97$, $\mu_1=0.5$, $\mu_2=0.3$, $\mu_3=0.2$, $b_1=108$
, $b_2=1.2$, $b_3=533/230$, $a_{11}=7$, $a_{12}=1$, $a_{13}=100$, $a_{21}=3/40$, $a_{22}=21/40$, $a_{23}=3/5$,
$a_{31}=27/230$, $a_{32}=1$, $a_{33}=6/5$,
system (\ref{seasonal-system}) satisfies the inequalities of the class 31 in Table \ref{biao0}. The numerical simulations for the solution of system \eqref{seasonal-system} with initial value $x_0=(0.7,1.6,0.8)$
and the orbit of the associated Poincar\'e map $\mathcal{P}$ are shown in Fig.  \ref{fig-c31}, which imply the system
admits an attracting invariant closed curve.}
\end{example}

\section{Discussion} \label{sec:7}
In this paper, we focus on the dynamics of the $3$-dimensional Lotka-Volterra competition system \eqref{seasonal-system} with seasonal succession. According to our research, system \eqref{seasonal-system} with seasonal succession is actually a very charming system,
which can be viewed as a bridge between the concrete discrete-time competitive mappings and time-periodically forced differential equations.
Our result seems to be the first attempt to study the classification of global dynamics for this special periodically forced differential equations with seasonal succession.

For system \eqref{seasonal-system},
based on the existence of a carrying simplex, we first propose an index formula on the sum of indices of all the
fixed points on the carrying simplex for the Poincar\'e map $\mathcal{P}$ associated with the system \eqref{seasonal-system}. The formula is similar to that provided in \cite{jiang2014} for competitive mappings, but our approach is totally different which
avoids constructing a competitive vector field. Generally, the method in \cite{jiang2014} is not applicable to the Poincar\'e map $\mathcal{P}$ because the explicit expression for  $\mathcal{P}$ is usually impossible to obtain, which makes the assumptions in \cite{jiang2014} difficult to verify for $\mathcal{P}$.
However, our method is easily applied to the competitive mappings with a carrying simplex.

By defining an equivalence relation relative to the local
dynamics of boundary fixed points, we derive the 33 stable equivalence classes of the dynamics for the Poincar\'e map $\mathcal{P}$.
 The parameter conditions and phase portrait for each class
are listed in Table \ref{biao0}. In classes 1--18, there is no positive fixed point and every trajectory tends to some fixed point
on the boundary. In classes 19--33, there is at least one (but not necessarily unique) positive fixed point. In particular, class 27 has a
heteroclinic cycle and we give the criteria on the stability of the heteroclinic cycle. Moreover, we prove that every orbit converges to some fixed point and obtain the
global dynamics for classes 19--25 and 33 if the positive fixed point is unique in these classes.  Our numerical experiments show that attracting invariant closed curves can occur in classes 26, 27, 29 and 31, on which all orbits are dense, that is, the associated solutions of system \eqref{seasonal-system} are quasiperiodic.

An interesting question that arise from our work need to be addressed. When or whether is the positive fixed point unique for the Poincar\'e map $\mathcal{P}$? We shall focus on this problem in the subsequent work \cite{nwx2023b}, where we provide conditions for the uniqueness of the positive fixed point and prove that there do exist some classes which can have multiple positive fixed points. This means that the uniqueness of the positive fixed point does not always hold for the Poincar\'e map  of system \eqref{seasonal-system}.

\section*{Acknowledgments}
The authors are very grateful to Prof. Sze-Bi Hsu for his valuable and useful discussions and suggestions.

\appendix
\section{The stable equivalence classes in $\mathrm{CLVS}(3)$}\label{appendix}
\setcounter{table}{0}
\begin{center}
  \begin{longtable}{c@{\extracolsep{\fill}}c@{\extracolsep{\fill}}c@{\extracolsep{\fill}}}
\caption{The $33$ equivalence classes in $\mathrm{CLVS}(3)$, where
$\gamma_{ij}=a_{ii}r_j-a_{ji}r_i$, $\beta_{ij}=\frac{a_{jj}r_i-a_{ij}r_j}{a_{ii}a_{jj}-a_{ij}a_{ji}}$, $i,j=1,2,3$, $i\neq j$. A fixed point is represented by a closed dot $\bullet$ if it attracts on the carrying simplex, by an open dot $\circ$ if it repels, and by the intersection of its stable and unstable manifolds if it is a saddle. The circle full of slashes in classes 19--33 denotes a region of unknown dynamics where there might be more than one fixed points or other complex dynamics such as invariant closed curves.} \\[-2pt]
        \hline
         \footnotesize {Classes} & \footnotesize {Parameter conditions} & \footnotesize {Phase Portraits} \\
        \hline
        \endfirsthead
        \caption[]{(continued)}\\
        \hline
        \footnotesize {Classes} & \footnotesize {Parameter conditions} & \footnotesize {Phase Portraits}   \\
        \hline
&&\\
        \endhead
        \hline
        \endfoot
        \endlastfoot
&&\\
1 &
\begin{tabular}{ll} &
{\footnotesize$\gamma_{12}<0, \gamma_{13}<0, \gamma_{21}>0$,}\\
& {\footnotesize$\gamma_{23}>0, \gamma_{31}>0, \gamma_{32}<0$}
\end{tabular}
&
\parbox{2cm}{\vspace{2pt}\includegraphics[width=2cm,height=1.6cm]{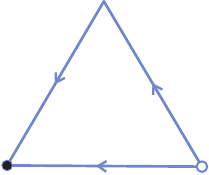}}
\\
&&\\[-10pt]
    2 &
\begin{tabular}{@{}l@{~}l@{}} {\footnotesize
 (i)}& {\footnotesize  $\gamma_{12}<0, \gamma_{13}<0, \gamma_{21}<0$,}\\
&{\footnotesize   $\gamma_{23}>0, \gamma_{31}>0, \gamma_{32}<0$}\\
{\footnotesize  (ii)}&{\footnotesize   $a_{31}\beta_{12}+a_{32}\beta_{21}<r_3$
} \end{tabular}
 &
    \parbox{2cm}{\vspace{2pt}\includegraphics[width=2cm,height=1.6cm]{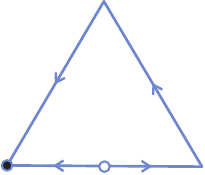}}\\
&&\\[-10pt]
    3 &
\begin{tabular}{ll} {\footnotesize
 (i)}&{\footnotesize  $\gamma_{12}<0, \gamma_{13}<0, \gamma_{21}>0$,}\\
 &{\footnotesize  $\gamma_{23}<0, \gamma_{31}>0, \gamma_{32}<0$}\\
{\footnotesize  (ii)}&{\footnotesize  $a_{12}\beta_{23}+a_{13}\beta_{32}<r_1$
} \end{tabular}
 &
    \parbox{2cm}{\vspace{2pt}\includegraphics[width=2cm,height=1.6cm]{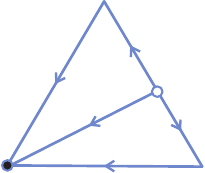}\vspace{2pt}}
\\
&&\\[-10pt]
4 &
\begin{tabular}{ll} {\footnotesize
 (i)}&{\footnotesize  $\gamma_{12}>0, \gamma_{13}<0, \gamma_{21}>0$,}\\
&{\footnotesize   $\gamma_{23}<0, \gamma_{31}>0, \gamma_{32}<0$}\\
{\footnotesize  (ii)}&{\footnotesize  $a_{12}\beta_{23}+a_{13}\beta_{32}<r_1$}\\
{\footnotesize  (iii)}&{\footnotesize  $a_{31}\beta_{12}+a_{32}\beta_{21}>r_3$
} \end{tabular}
 &
    \parbox{2cm}{\vspace{2pt}\includegraphics[width=2cm,height=1.6cm]{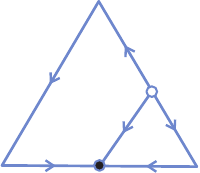}\vspace{2pt}}\\
&&\\[-10pt]
    5 &
\begin{tabular}{ll} {\footnotesize
 (i)}&{\footnotesize  $\gamma_{12}>0, \gamma_{13}>0, \gamma_{21}>0$,}\\
&{\footnotesize   $\gamma_{23}<0, \gamma_{31}<0, \gamma_{32}>0$}\\
{\footnotesize  (ii)}&{\footnotesize  $a_{31}\beta_{12}+a_{32}\beta_{21}>r_3$
} \end{tabular}
 &
    \parbox{2cm}{\vspace{2pt}\includegraphics[width=2cm,height=1.6cm]{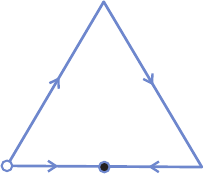}\vspace{2pt}}\\
    &&\\[-10pt]
6 &
\begin{tabular}{ll} {\footnotesize
 (i)}&{\footnotesize  $\gamma_{12}>0, \gamma_{13}>0, \gamma_{21}<0$,}\\
&{\footnotesize   $\gamma_{23}>0, \gamma_{31}<0, \gamma_{32}>0$}\\
{\footnotesize  (ii)}&{\footnotesize  $a_{12}\beta_{23}+a_{13}\beta_{32}>r_1$
} \end{tabular}
 &
\parbox{2cm}{\vspace{2pt}\includegraphics[width=2cm,height=1.6cm]{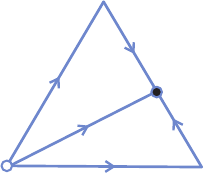}\vspace{2pt}}
\\
&&\\[-10pt]
    7 &
\begin{tabular}{ll} {\footnotesize
 (i)}&{\footnotesize  $\gamma_{12}>0, \gamma_{13}>0, \gamma_{21}>0$,}\\
&{\footnotesize   $\gamma_{23}>0, \gamma_{31}<0, \gamma_{32}<0$}\\
{\footnotesize  (ii)}&{\footnotesize  $a_{31}\beta_{12}+a_{32}\beta_{21}<r_3$}
 \end{tabular}
 &
\parbox{2cm}{\vspace{2pt}\includegraphics[width=2cm,height=1.6cm]{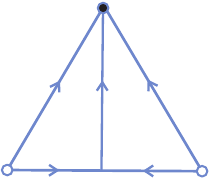}\vspace{2pt}}
\\
&&\\[-10pt]
8 &
\begin{tabular}{ll} {\footnotesize
 (i)}&{\footnotesize  $\gamma_{12}>0, \gamma_{13}>0, \gamma_{21}>0$,}\\
&{\footnotesize   $\gamma_{23}<0, \gamma_{31}<0, \gamma_{32}<0$}\\
{\footnotesize  (ii)}&{\footnotesize  $a_{12}\beta_{23}+a_{13}\beta_{32}<r_1$}\\
{\footnotesize  (iii)}&{\footnotesize  $a_{31}\beta_{12}+a_{32}\beta_{21}<r_3$
} \end{tabular}
 &
\parbox{2cm}{\vspace{2pt}\includegraphics[width=2cm,height=1.6cm]{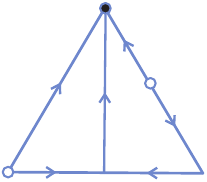}\vspace{2pt}}
\\
&&\\[-10pt]
    9 &
\begin{tabular}{ll} {\footnotesize
 (i)}&{\footnotesize  $\gamma_{12}>0, \gamma_{13}>0, \gamma_{21}>0$,}\\
&{\footnotesize   $\gamma_{23}>0, \gamma_{31}<0, \gamma_{32}>0$}\\
{\footnotesize  (ii)}&{\footnotesize  $a_{12}\beta_{23}+a_{13}\beta_{32}>r_1$}\\
{\footnotesize  (iii)}&{\footnotesize  $a_{31}\beta_{12}+a_{32}\beta_{21}<r_3$
} \end{tabular}
 &
\parbox{2cm}{\vspace{2pt}\includegraphics[width=2cm,height=1.6cm]{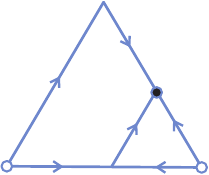}\vspace{2pt}}
\\
&&\\[-10pt]
10 &
\begin{tabular}{ll} {\footnotesize
 (i)}&{\footnotesize  $\gamma_{12}>0, \gamma_{13}>0, \gamma_{21}>0$,}\\
&{\footnotesize   $\gamma_{23}>0, \gamma_{31}<0, \gamma_{32}>0$}\\
{\footnotesize  (ii)}&{\footnotesize  $a_{12}\beta_{23}+a_{13}\beta_{32}<r_1$}\\
{\footnotesize  (iii)}&{\footnotesize  $a_{31}\beta_{12}+a_{32}\beta_{21}>r_3$
} \end{tabular}
 &
\parbox{2cm}{\vspace{2pt}\includegraphics[width=2cm,height=1.6cm]{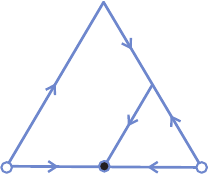}\vspace{2pt}}
\\
&&\\[-10pt]
    11 &
\begin{tabular}{ll} {\footnotesize
 (i)}&{\footnotesize  $\gamma_{12}>0, \gamma_{13}>0, \gamma_{21}>0$,}\\
&{\footnotesize   $\gamma_{23}<0, \gamma_{31}>0, \gamma_{32}<0$}\\
{\footnotesize  (ii)}& {\footnotesize  $a_{12}\beta_{23}+a_{13}\beta_{32}<r_1$}\\
{\footnotesize  (iii)}& {\footnotesize  $a_{21}\beta_{13}+a_{23}\beta_{31}<r_2$}\\
{\footnotesize  (iv)}& {\footnotesize  $a_{31}\beta_{12}+a_{32}\beta_{21}>r_3$
} \end{tabular}
 &
\parbox{2cm}{\vspace{2pt}\includegraphics[width=2cm,height=1.6cm]{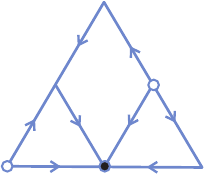}\vspace{2pt}}
\\
&&\\[-10pt]
12 &
\begin{tabular}{ll} {\footnotesize
 (i)}&{\footnotesize  $\gamma_{12}>0, \gamma_{13}>0, \gamma_{21}>0$,}\\
&{\footnotesize   $\gamma_{23}>0, \gamma_{31}>0, \gamma_{32}>0$}\\
{\footnotesize  (ii)}&{\footnotesize  $a_{12}\beta_{23}+a_{13}\beta_{32}<r_1$}\\
{\footnotesize  (iii)}&{\footnotesize  $a_{21}\beta_{13}+a_{23}\beta_{31}<r_2$}\\
{\footnotesize  (iv)}&{\footnotesize  $a_{31}\beta_{12}+a_{32}\beta_{21}>r_3$
} \end{tabular}
 &
\parbox{2cm}{\vspace{2pt}\includegraphics[width=2cm,height=1.6cm]{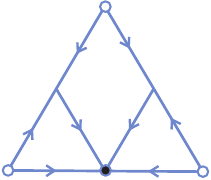}\vspace{2pt}}
\\
&&\\[-10pt]
    13 &
\begin{tabular}{ll} {\footnotesize
 (i)}&{\footnotesize  $\gamma_{12}<0, \gamma_{13}<0, \gamma_{21}<0$,}\\
&{\footnotesize   $\gamma_{23}<0, \gamma_{31}>0, \gamma_{32}>0$}\\
{\footnotesize  (ii)}&{\footnotesize  $a_{31}\beta_{12}+a_{32}\beta_{21}>r_3$
} \end{tabular}
 &
\parbox{2cm}{\vspace{2pt}\includegraphics[width=2cm,height=1.6cm]{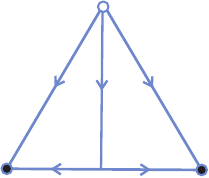}\vspace{2pt}}
\\
&&\\[-10pt]
14 &
\begin{tabular}{ll} {\footnotesize
 (i)}&{\footnotesize  $\gamma_{12}<0, \gamma_{13}<0, \gamma_{21}<0$,}\\
&{\footnotesize   $\gamma_{23}>0, \gamma_{31}>0, \gamma_{32}>0$}\\
{\footnotesize  (ii)}&{\footnotesize  $a_{12}\beta_{23}+a_{13}\beta_{32}>r_1$}\\
{\footnotesize  (iii)}&{\footnotesize  $a_{31}\beta_{12}+a_{32}\beta_{21}>r_3$
} \end{tabular}
 &
\parbox{2cm}{\vspace{2pt}\includegraphics[width=2cm,height=1.6cm]{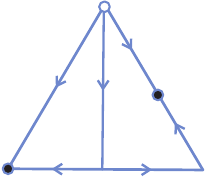}\vspace{2pt}}
\\
&&\\[-10pt]
    15 &
\begin{tabular}{ll} {\footnotesize
 (i)}&{\footnotesize  $\gamma_{12}<0, \gamma_{13}<0, \gamma_{21}<0$,}\\
 &{\footnotesize  $\gamma_{23}<0, \gamma_{31}>0, \gamma_{32}<0$}\\
{\footnotesize  (ii)}&{\footnotesize  $a_{12}\beta_{23}+a_{13}\beta_{32}<r_1$}\\
{\footnotesize  (iii)}&{\footnotesize  $a_{31}\beta_{12}+a_{32}\beta_{21}>r_3$
} \end{tabular}
 &
\parbox{2cm}{\vspace{2pt}\includegraphics[width=2cm,height=1.6cm]{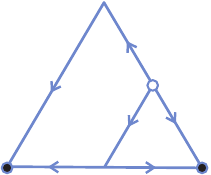}\vspace{2pt}}
 \\
&&\\[-10pt]
16 &
\begin{tabular}{ll} {\footnotesize
 (i)}&{\footnotesize  $\gamma_{12}<0, \gamma_{13}<0, \gamma_{21}<0$,}\\
&{\footnotesize   $\gamma_{23}<0, \gamma_{31}>0, \gamma_{32}<0$}\\
{\footnotesize  (ii)}&{\footnotesize  $a_{12}\beta_{23}+a_{13}\beta_{32}>r_1$}\\
{\footnotesize  (iii)}&{\footnotesize  $a_{31}\beta_{12}+a_{32}\beta_{21}<r_3$
} \end{tabular}
 &
\parbox{2cm}{\vspace{2pt}\includegraphics[width=2cm,height=1.6cm]{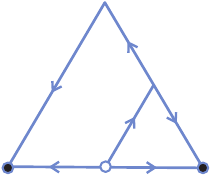}\vspace{2pt}} \\
&&\\[-10pt]
    17 &
\begin{tabular}{ll} {\footnotesize
 (i)}&{\footnotesize  $\gamma_{12}<0, \gamma_{13}<0, \gamma_{21}<0$,}\\
&{\footnotesize   $\gamma_{23}>0, \gamma_{31}<0, \gamma_{32}>0$}\\
{\footnotesize  (ii)}&{\footnotesize  $a_{12}\beta_{23}+a_{13}\beta_{32}>r_1$}\\
{\footnotesize  (iii)}&{\footnotesize  $a_{21}\beta_{13}+a_{23}\beta_{31}>r_2$}\\
{\footnotesize  (iv)}&{\footnotesize  $a_{31}\beta_{12}+a_{32}\beta_{21}<r_3$
} \end{tabular}
 &
\parbox{2cm}{\vspace{2pt}\includegraphics[width=2cm,height=1.6cm]{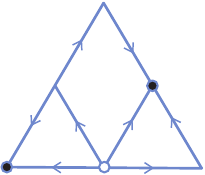}\vspace{2pt}}
 \\
&&\\[-10pt]
    18 &
\begin{tabular}{ll} {\footnotesize
 (i)}&{\footnotesize  $\gamma_{12}<0, \gamma_{13}<0, \gamma_{21}<0$,}\\
&{\footnotesize   $\gamma_{23}<0, \gamma_{31}<0, \gamma_{32}<0$}\\
{\footnotesize  (ii)}&{\footnotesize  $a_{12}\beta_{23}+a_{13}\beta_{32}>r_1$}\\
{\footnotesize  (iii)}& {\footnotesize  $a_{21}\beta_{13}+a_{23}\beta_{31}>r_2$}\\
{\footnotesize  (iv)}& {\footnotesize  $a_{31}\beta_{12}+a_{32}\beta_{21}<r_3$
} \end{tabular}
 &
\parbox{2cm}{\vspace{2pt}\includegraphics[width=2cm,height=1.6cm]{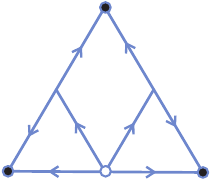}\vspace{2pt}}
 \\
&&\\[-10pt]
19 &
\begin{tabular}{ll} {\footnotesize
 (i)}&{\footnotesize  $\gamma_{12}>0, \gamma_{13}>0, \gamma_{21}<0$,}\\
&{\footnotesize   $\gamma_{23}<0, \gamma_{31}<0, \gamma_{32}<0$}\\
{\footnotesize  (ii)}&{\footnotesize  $a_{12}\beta_{23}+a_{13}\beta_{32}<r_1$
} \end{tabular}
 &
\parbox{2cm}{\vspace{2pt}\includegraphics[width=2cm,height=1.6cm]{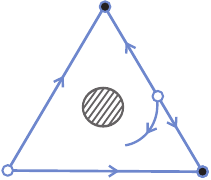}\vspace{2pt}}
 \\
&&\\[-10pt]
    20 &
\begin{tabular}{ll} {\footnotesize
 (i)}&{\footnotesize  $\gamma_{12}<0, \gamma_{13}<0, \gamma_{21}<0$,}\\
&{\footnotesize   $\gamma_{23}<0, \gamma_{31}>0, \gamma_{32}<0$}\\
{\footnotesize  (ii)}&{\footnotesize  $a_{12}\beta_{23}+a_{13}\beta_{32}<r_1$}\\
{\footnotesize  (iii)}&{\footnotesize  $a_{31}\beta_{12}+a_{32}\beta_{21}<r_3$
} \end{tabular}
 &
\parbox{2cm}{\vspace{2pt}\includegraphics[width=2cm,height=1.6cm]{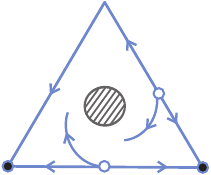}\vspace{2pt}}
\\
&&\\[-10pt]
21 &
\begin{tabular}{ll} {\footnotesize
 (i)}&{\footnotesize  $\gamma_{12}<0, \gamma_{13}<0, \gamma_{21}<0$,}\\
&{\footnotesize   $\gamma_{23}>0, \gamma_{31}<0, \gamma_{32}>0$}\\
{\footnotesize  (ii)}&{\footnotesize  $a_{12}\beta_{23}+a_{13}\beta_{32}>r_1$}\\
{\footnotesize  (iii)}&{\footnotesize  $a_{21}\beta_{13}+a_{23}\beta_{31}<r_2$}\\
{\footnotesize  (iv)}& {\footnotesize  $a_{31}\beta_{12}+a_{32}\beta_{21}<r_3$
} \end{tabular}
 &
\parbox{2cm}{\vspace{2pt}\includegraphics[width=2cm,height=1.6cm]{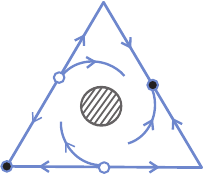}\vspace{2pt}}
\\
&&\\[-10pt]
    22 &
\begin{tabular}{ll} {\footnotesize
 (i)}&{\footnotesize  $\gamma_{12}>0, \gamma_{13}>0, \gamma_{21}<0$,}\\
&{\footnotesize   $\gamma_{23}<0, \gamma_{31}>0, \gamma_{32}<0$}\\
{\footnotesize  (ii)}&{\footnotesize  $a_{12}\beta_{23}+a_{13}\beta_{32}<r_1$}\\
{\footnotesize  (iii)}& {\footnotesize  $a_{21}\beta_{13}+a_{23}\beta_{31}>r_2$
} \end{tabular}
 &
\parbox{2cm}{\vspace{2pt}\includegraphics[width=2cm,height=1.6cm]{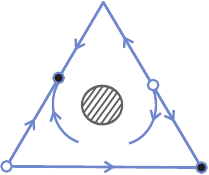}\vspace{2pt}}
\\
&&\\[-10pt]
23 &
\begin{tabular}{ll} {\footnotesize
 (i)}&{\footnotesize  $\gamma_{12}>0, \gamma_{13}>0, \gamma_{21}>0$,}\\
&{\footnotesize   $\gamma_{23}>0, \gamma_{31}<0, \gamma_{32}<0$}\\
{\footnotesize  (ii)}&{\footnotesize  $a_{31}\beta_{12}+a_{32}\beta_{21}>r_3$
} \end{tabular}
 &
\parbox{2cm}{\vspace{2pt}\includegraphics[width=2cm,height=1.6cm]{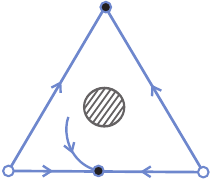}\vspace{2pt}}
\\
&&\\[-10pt]
    24 &
\begin{tabular}{ll} {\footnotesize
 (i)}&{\footnotesize  $\gamma_{12}>0, \gamma_{13}>0, \gamma_{21}>0$,}\\
&{\footnotesize   $\gamma_{23}>0, \gamma_{31}<0, \gamma_{32}>0$}\\
{\footnotesize  (ii)}&{\footnotesize  $a_{12}\beta_{23}+a_{13}\beta_{32}>r_1$}\\
{\footnotesize  (iii)}&{\footnotesize  $a_{31}\beta_{12}+a_{32}\beta_{21}>r_3$
} \end{tabular}
&
\parbox{2cm}{\vspace{2pt}\includegraphics[width=2cm,height=1.6cm]{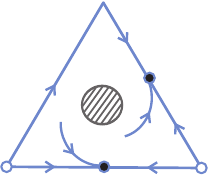}\vspace{2pt}}
 \\
&&\\[-10pt]
25 &
\begin{tabular}{ll} {\footnotesize
 (i)}&{\footnotesize  $\gamma_{12}>0, \gamma_{13}>0, \gamma_{21}>0$,}\\
&{\footnotesize   $\gamma_{23}<0, \gamma_{31}>0, \gamma_{32}<0$}\\
{\footnotesize  (ii)}&{\footnotesize  $a_{12}\beta_{23}+a_{13}\beta_{32}<r_1$}\\
{\footnotesize  (iii)}&{\footnotesize  $a_{21}\beta_{13}+a_{23}\beta_{31}>r_2$}\\
{\footnotesize  (iv)}&{\footnotesize  $a_{31}\beta_{12}+a_{32}\beta_{21}>r_3$
} \end{tabular}
 &
\parbox{2cm}{\vspace{2pt}\includegraphics[width=2cm,height=1.6cm]{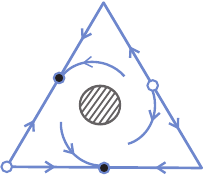}\vspace{2pt}}
\\
&&\\[-10pt]
    26 &
\begin{tabular}{ll} {\footnotesize
 (i)}&{\footnotesize  $\gamma_{12}>0, \gamma_{13}>0, \gamma_{21}<0$,}\\
&{\footnotesize   $\gamma_{23}<0, \gamma_{31}>0, \gamma_{32}<0$}\\
{\footnotesize  (ii)}&{\footnotesize  $a_{12}\beta_{23}+a_{13}\beta_{32}>r_1$}\\
{\footnotesize  (iii)}&{\footnotesize  $a_{21}\beta_{13}+a_{23}\beta_{31}<r_2$
} \end{tabular}
 &
\parbox{2cm}{\vspace{2pt}\includegraphics[width=2cm,height=1.6cm]{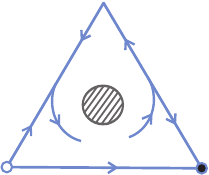}\vspace{2pt}}\\
&&\\[-10pt]
27 &
\begin{tabular}{ll} &{\footnotesize
 $\gamma_{12}>0, \gamma_{13}<0, \gamma_{21}<0$,}\\
&{\footnotesize   $\gamma_{23}>0, \gamma_{31}>0, \gamma_{32}<0$
} \end{tabular}
 &
    \parbox{2cm}{\vspace{2pt}\includegraphics[width=2cm,height=1.6cm]{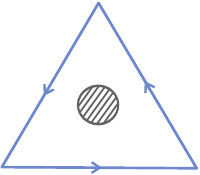}\vspace{2pt}}\\
&&\\[-10pt]
    28 &
\begin{tabular}{ll} {\footnotesize
 (i)}&{\footnotesize  $\gamma_{12}<0, \gamma_{13}<0, \gamma_{21}<0$,}\\
 &{\footnotesize  $\gamma_{23}>0, \gamma_{31}>0, \gamma_{32}<0$}\\
{\footnotesize  (ii)}&{\footnotesize   $a_{31}\beta_{12}+a_{32}\beta_{21}>r_3$
} \end{tabular}
&
\parbox{2cm}{\vspace{2pt}\includegraphics[width=2cm,height=1.6cm]{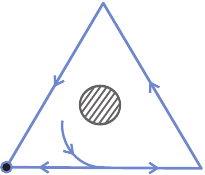}\vspace{2pt}}\\
&&\\[-10pt]
29 &
\begin{tabular}{ll} {\footnotesize
 (i)}&{\footnotesize  $\gamma_{12}>0, \gamma_{13}>0, \gamma_{21}>0$,}\\
&{\footnotesize   $\gamma_{23}<0, \gamma_{31}<0, \gamma_{32}>0$}\\
{\footnotesize  (ii)}&{\footnotesize  $a_{31}\beta_{12}+a_{32}\beta_{21}<r_3$
} \end{tabular}
 &
\parbox{2cm}{\vspace{2pt}\includegraphics[width=2cm,height=1.6cm]{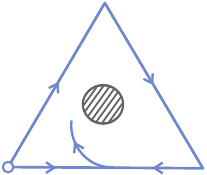}\vspace{2pt}}\\
&&\\[-10pt]
    30 &
\begin{tabular}{ll} {\footnotesize
(i)}&{\footnotesize  $\gamma_{12}<0, \gamma_{13}<0, \gamma_{21}<0$,}\\
&{\footnotesize  $\gamma_{23}<0, \gamma_{31}>0, \gamma_{32}<0$}\\
{\footnotesize  (ii)}&{\footnotesize  $a_{12}\beta_{23}+a_{13}\beta_{32}>r_1$}\\
{\footnotesize  (iii)}&{\footnotesize  $a_{31}\beta_{12}+a_{32}\beta_{21}>r_3$
} \end{tabular}
&
\parbox{2cm}{\vspace{2pt}\includegraphics[width=2cm,height=1.6cm]{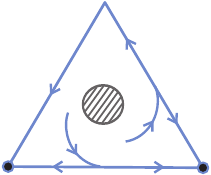}\vspace{2pt}} \\
&&\\[-10pt]
    31 &
\begin{tabular}{ll} {\footnotesize
(i)}&{\footnotesize  $\gamma_{12}>0, \gamma_{13}>0, \gamma_{21}>0$,}\\
&{\footnotesize  $\gamma_{23}>0, \gamma_{31}<0, \gamma_{32}>0$}\\
{\footnotesize  (ii)}&{\footnotesize  $a_{12}\beta_{23}+a_{13}\beta_{32}<r_1$}\\
{\footnotesize  (iii)}&{\footnotesize  $a_{31}\beta_{12}+a_{32}\beta_{21}<r_3$
} \end{tabular}
&
\parbox{2cm}{\vspace{2pt}\includegraphics[width=2cm,height=1.6cm]{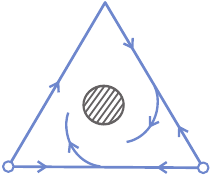}\vspace{2pt}}\\
&&\\[-10pt]
32 &
\begin{tabular}{ll} {\footnotesize
 (i)}&{\footnotesize  $\gamma_{12}<0, \gamma_{13}<0, \gamma_{21}<0$,}\\
&{\footnotesize   $\gamma_{23}<0, \gamma_{31}<0, \gamma_{32}<0$}\\
{\footnotesize  (ii)}&{\footnotesize  $a_{12}\beta_{23}+a_{13}\beta_{32}>r_1$}\\
{\footnotesize  (iii)}&{\footnotesize  $a_{21}\beta_{13}+a_{23}\beta_{31}>r_2$}\\
{\footnotesize  (iv)}&{\footnotesize  $a_{31}\beta_{12}+a_{32}\beta_{21}>r_3$
} \end{tabular}
 &
\parbox{2cm}{\vspace{2pt}\includegraphics[width=2cm,height=1.6cm]{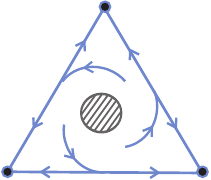}\vspace{2pt}}\\
&&\\[-10pt]
    33 &
\begin{tabular}{ll} {\footnotesize
 (i)}&{\footnotesize  $\gamma_{12}>0, \gamma_{13}>0, \gamma_{21}>0$,}\\
&{\footnotesize   $\gamma_{23}>0, \gamma_{31}>0, \gamma_{32}>0$}\\
{\footnotesize  (ii)}&{\footnotesize  $a_{12}\beta_{23}+a_{13}\beta_{32}<r_1$}\\
{\footnotesize  (iii)}&{\footnotesize  $a_{21}\beta_{13}+a_{23}\beta_{31}<r_2$}\\
{\footnotesize  (iv)}&{\footnotesize  $a_{31}\beta_{12}+a_{32}\beta_{21}<r_3$
} \end{tabular}
&
\parbox{2cm}{\vspace{2pt}\includegraphics[width=2cm,height=1.6cm]{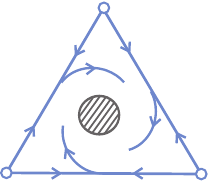}\vspace{2pt}}\\
[-8pt]\label{biao0}
\end{longtable}
\end{center}
\section{The phase portraits for classes 19--25 and 33 when the positive fixed point is unique}\label{appendix-B}

\begin{figure}[h]
	\begin{center}
		\includegraphics[width=0.66\textwidth]{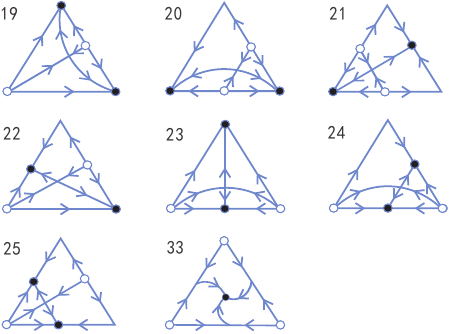}
		\caption{The phase portraits on the carrying simplices for classes 19--25 and 33 whenever there is a unique positive fixed point for these classes. The fixed point notation is as in Table \ref{biao0}.} \label{fig:c19-25-33}
	\end{center}
\end{figure}

\bibliographystyle{siamplain}
\bibliography{refs}
\end{document}